\documentclass{amsart}
\usepackage{amssymb}
\usepackage{amsthm}
\usepackage{graphics}
\newtheorem*{reftheorem}{Theorem}
\newtheorem{theorem}{Theorem}[section]
\newtheorem{lemma}[theorem]{Lemma}

\newtheorem*{refcor}{Corollary}

\numberwithin{equation}{section}
\title{Lipshitz Maps from Surfaces}
\author{Larry Guth}
\address{Department of Mathematics, M.I.T., Cambridge MA 02139, USA}
\email{larryg@math.mit.edu}

\begin{document}
\begin{abstract}
We give a simple procedure to estimate the smallest Lipshitz
constant of a degree 1 map from a Riemannian 2-sphere to the unit
2-sphere, up to a factor of 10.  Using this procedure, we are
able to prove several inequalities involving this Lipshitz
constant.  For instance, if the smallest Lipshitz constant is at
least 1, then the Riemannian 2-sphere has Uryson 1-Width less
than 12 and contains a closed geodesic of length less than 160. 
Similarly, if a closed oriented Riemannian surface does not
admit a degree 1 map to the unit 2-sphere with Lipshitz constant
1, then it contains a closed homologically non-trivial curve of
length less than $4 \pi$.  On the other hand, we give examples of
high genus surfaces with arbitrarily large Uryson 1-Width which
do not admit a map of non-zero degree to the unit sphere with
Lipshitz constant 1. 
\end{abstract}

\maketitle

This paper is about estimating the best Lipshitz constant of a
degree 1 map from a closed oriented Riemannian surface to the
unit 2-sphere.  The difficulty of this problem depends on the genus
of the surface.  In case the surface is a topological 2-sphere,
we give a simple procedure for estimating the best Lipshitz
constant within a factor of twelve.

Pick a point p on the Riemannian 2-sphere $(S^2, g)$.  Consider
the distance spheres around p.  Let D be the largest diameter of
any connected component of any of the distance spheres around p.

\begin{theorem} The best Lipshitz constant of a degree 1 map from
$(S^2, g)$ to the unit 2-sphere is more than $1/D$ and less than
$12 / D$.
\end{theorem}

We now recall some vocabulary, which we will use all through the
paper, designed to describe how large and how wide Riemannian
manifolds are.  The hypersphericity of a Riemannian n-manifold M
is the supremal R so that there is a contracting map of non-zero
degree from M to the n-sphere of radius R.  In this paper, we
will mostly focus on degree 1 maps, so we define the degree 1
hypersphericity of a Riemannian n-manifold M as the supremal R so
that there is a contracting map of degree 1 from M to the
n-sphere of radius R.  The Uryson k-Width of a metric space X is
the infimal W so that there is a continuous map from X to a
k-dimensional polyhedron whose fibers have diameter less than W. 
The main theorem about hypersphericity and Uryson Width is an
estimate by Gromov in \cite{G1} that the hypersphericity of an
n-manifold is less than its Uryson (n-1)-width.

Rephrased in this vocabulary, our first theorem states that the
degree 1 hypersphericity of $(S^2, g)$ is between D/12 and D. 
Using Gromov's result, it follows that the Uryson 1-Width of
$(S^2, g)$ is also between D/12 and D.  In particular, we see
that the hypersphericity and the Uryson 1-Width of a 2-sphere
agree up to a factor of twelve.  Since it is not hard to write an
efficient algorithm to estimate D to any desired accuracy, we can
efficiently estimate the hypersphericity and Uryson 1-width of a
Riemannian 2-sphere up to a factor of twelve.

We prove Theorem 0.1 in section 1.  In the next three
sections, we apply the first theorem to give estimates relating
the degree 1 hypersphericity of a surface to other geometric
quantities.  While the first theorem is quite easy, some of these
applications are more difficult.  In the last section of the
paper, we show that the first theorem fails for surfaces of high
genus. 

In section 2, we bound the length of the shortest closed geodesic
on a 2-sphere by its degree 1 hypersphericity.

\begin{theorem} If $(S^2, g)$ has degree 1 hypersphericity less
than 1, then it contains a closed geodesic of length
less than 160.
\end{theorem}

This theorem extends a result of Croke, who proved in \cite{C}
that if the area of a sphere is less than $4 \pi$, then it
contains a closed geodesic of length less than 100.  The main
step in the proof is to produce a sufficiently nice family of
short curves on $(S^2, g)$.  In particular, we prove along the
way the following theorem.

\begin{theorem} If $(S^2, g)$ has degree 1 hypersphericity
less than 1, then it admits a map to a tree whose fibers have
length less than 120.
\end{theorem}

By our first theorem, we already know that $(S^2, g)$ admits a
map to a tree whose fibers have diameter less than 12.  It may
happen, however, that the fibers are extremely long curves with
many wiggles.  Roughly speaking, the above theorem means that it
is possible to straighten out all of these wiggles (possibly
after replacing the target tree by a tree with much more
branching).

In section 3, we apply our estimates to surfaces in Euclidean
space.  Our first estimate says that if we slice the unit ball
into non-singular surfaces, then one of these surfaces has a
large hypersphericity.

\begin{theorem} Let f be an embedding of the unit n-ball into
$\mathbb{R}^n$, and let $f_0$ be the first n-2 coordinates of f. 
The level sets of $f_0$ are planar surfaces in the unit n-ball. 
One of the level sets of $f_0$ has a connected component with
degree 1 hypersphericity at least 1/24.
\end{theorem}

Our second estimate is a variation of the isoperimetric
inequality.  For a 2-sphere embedded in $\mathbb{R}^3$, our
inequality controls the volume inside the 2-sphere in terms of
the hypersphericities of pieces of the boundary.

\begin{theorem} Let U be a bounded open set in $\mathbb{R}^3$, with
boundary diffeomorphic to a 2-sphere.  Then the boundary of U
contains disjoint subsets $S_i$, so that the following inequality
holds.

$$\rm{Vol}(U) < 10^6 \sum \rm{HS}(S_i)^3.$$

In this equation, $\rm{HS}(S_i)$ stands for the degree 1
hypersphericity of $S_i$.
\end{theorem}

Our third estimate gives a lower bound for the hypersphericity of
the boundary of U.

\begin{theorem} Let U be a bounded open set in $\mathbb{R}^3$ with
boundary diffeomorphic to a 2-sphere, and suppose that U admits
an area-contracting diffeomorphism to the unit ball.  Then the
degree 1 hypersphericity of the boundary of U is at least
1/200.
\end{theorem}

The constants in the previous two theorems are probably rather
poor.  In particular, it seems possible that if U admits an
area-contracting diffeomorphism to the unit ball, then the
hypersphericity of the boundary of U is at least 1.

Section 4 proves some estimates for surfaces of higher genus.  We
first prove two theorems bounding the lengths of homologically
interesting curves in terms of the hypersphericity.

\begin{theorem} Let $(\Sigma,g)$ be a closed, oriented surface of
genus at least 1 and degree 1 hypersphericity less than 1.  Then
it contains a homologically non-trivial closed curve of length
less than $4 \pi$.
\end{theorem}

\begin{theorem} Let $(\Sigma, g)$ be a closed oriented surface of
genus G and degree 1 hypersphericity less than 1.  Then it
contains homologically independent curves $C_1$, ..., $C_G$ with
the length of $C_k$ bounded by $200k$. \end{theorem}

By a similar method, we bound the Uryson 1-Width of a surface
in terms of its degree 1 hypersphericity and its genus.

\begin{theorem} Let $(\Sigma, g)$ be a closed oriented surface of
genus G and degree 1 hypersphericity less than 1.  Then the
Uryson 1-Width of $(\Sigma, g)$ is less than $200 G + 12$.
\end{theorem}

The power of G in this estimate may not be sharp, but in
section 5 we will construct examples of surfaces with degree 1
hypersphericity 1 and Uryson 1-Width on the order of $G^{1/2}$,
and surfaces with hypersphericity 1 and Uryson 1-Width on the
order of $G^{1/4}$.

In spite of the results in section 4, I don't really understand
what a high genus surface with small hypersphericity is like. 
For example, I am very far from being able to write an algorithm
that estimates the hypersphericity of a surface within a constant
factor.  In section 5, we will give examples of metrics on
surfaces of very high genus with arbitrarily small
hypersphericity and Uryson 1-Width at least 1.

\begin{theorem} There exist closed oriented Riemannian surfaces
with arbitrarily small hypersphericity and Uryson 1-Width at
least 1. \end{theorem}

We will give three different examples of surfaces proving this
theorem.  The three different examples have small hypersphericity
for three different ``reasons'', and I include them to suggest
that the problem of estimating the best Lipshitz constant of a
degree 1 map from a surface to the unit sphere is a complicated
problem, which can involve several different pieces of geometry
and topology.  Our three examples will use a little homotopy
theory, a little lattice theory, and a little bit of the theory
of systoles.

I would like to thank my thesis advisor, Tom Mrowka, for his
patience and support, and Mikhail Katz, who read an earlier version
of the paper and made helpful suggestions.  I am also grateful to
the referee for her or his careful reading and cogent suggestions.

\tableofcontents

\section{Hypersphericity of 2-Spheres}

In this section, we estimate the best Lipshitz constant of a
degree 1 map from a Riemannian 2-sphere $(S^2,g)$ to the unit
2-sphere.

Pick any point p in $(S^2, g)$.  Let $S(p, R)$ be the metric sphere
around p of radius R - that is, the set of points x in $(S^2, g)$
whose distance from p is R.  Let D be the the largest diameter of
any connected component of the distance sphere $S(p, R)$ for any
radius R.  We will prove that the smallest Lipshitz constant of any
degree 1 map from $(S^2, g)$ to the unit sphere is approximately
$5/D$.

\begin{theorem} The smallest Lipshitz constant of a degree 1 map from
$(S^2, g)$ to the unit 2-sphere is at least $\pi / (2D)$ and at
most $(2 + \sqrt2) \pi / D$.
\end{theorem}

\proof To prove the upper bound, we will construct a Lipshitz
degree 1 map.  For some radius R, one component of $S(p, R)$ has
diameter D.  Choose two points q and r in that component, so that
the distance from q to r is D.  We define $F_0(x) = (dist(p,x),
dist(q,x))$.  The function $F_0$ maps $(S^2, g)$ to $\mathbb{R}^2$, with
Lipshitz constant $\sqrt{2}$. 

We will consider the restriction of $F_0$ to three curves in
$S^2$: a minimal geodesic $g_{pq}$ from p to q, a minimal
geodesic $g_{pr}$ from p to r, and a curve $\gamma_{qr}$ joining
q to r.  For now, we will assume that the component of $S(p,R)$
containing q and r is path connected, and we take $\gamma_{qr}$
to lie in this component.

The image under $F_0$ of the geodesic $g_{pq}$ is the straight
line from $(0,R)$ to $(R,0)$.  Because $\gamma_{qr}$ lies in
$S(p,R)$, the image of $\gamma_{qr}$ under $F_0$ is contained in
the line $x = R$.  The image of the geodesic $g_{pr}$ cannot be
determined exactly, but we will prove that it lies above a
certain line.  Let $g_{pr}(t)$ parametrize $g_{pr}$ by arclength
with $g_{pr}(0) = p$ and $g_{pr}(R)=r$.  By definition, $dist(p,
g_{pr}(t))=t$.  Then $F_0 \circ g_{pr}(t) = (t, dist(q,
g_{pr}(t)))$, and the image $F_0(g_{pr})$ is the graph of the
function $f(t) = dist(q, g_{pr}(t))$.  The function f has
Lipshitz norm less than or equal to 1.  Since the distance from q
to r is D, $f(R) = D$.  Therefore, the graph of f lies above the
line from $(R - D/2, D/2)$ to $(R,D)$. 

\includegraphics{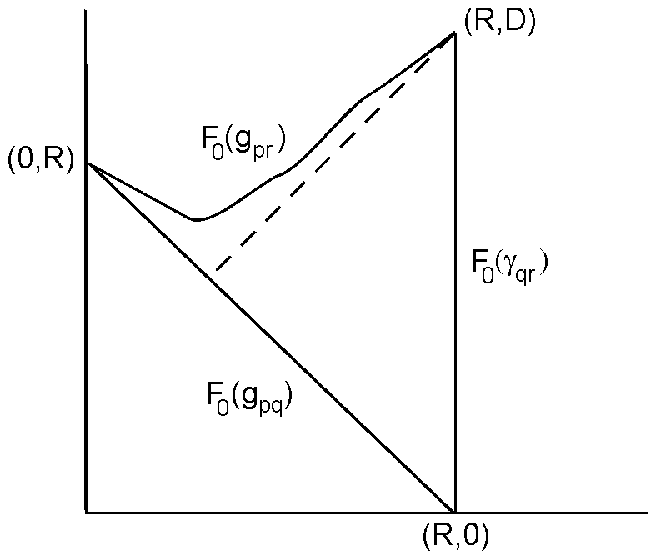}

Let T be the open triangle in $\mathbb{R}^2$ with vertices
$(R,0), (R,D),$ and $(R-D/2, D/2)$.  We have proved that
$F_0(g_{pr})$ does not intersect T, and it is easy to check that
$F_0(g_{pq})$ and $F_0(\gamma_{qr})$ don't intersect T either. 
Let c be the embedded circle in $(S^2, g)$ formed by joining
$g_{pq}$, $\gamma_{qr}$, and $g_{pr}$.  After choosing an
appropriate orientation of c, $F_0: c \rightarrow \mathbb{R}^2 -
T$ has degree one.

The triangle T contains an embedded disk U of radius $\frac{1}{2
\sqrt2 + 2} D$.  Using the exponential map, we construct a
homeomorphism $\pi_1$ from the disk U to the Northern hemisphere
of the unit sphere.  We define a map $\pi_2$ from the disk U to
the Southern hemisphere of the unit sphere by composing $\pi_1$
with a reflection through the plane of the equator.  We extend
$\pi_1$ and $\pi_2$ to all of $\mathbb{R}^2$ as follows.  For any
point x outside of U, let p(x) be the nearest point to x on the
boundary of U.  Now, define $\pi_1(x) = \pi_1(p(x))$ and
$\pi_2(x) = \pi_2(p(x))$.  The maps $\pi_1$ and $\pi_2$ agree on
the the boundary of U and hence on all the points outside of U. 
The maps $\pi_1$ and $\pi_2$ have Lipshitz constant $(\sqrt 2 +
1) \pi / D$.

Finally we define our degree 1 map.  The curve c divides $(S^2,
g)$ into two disks, $D_1$ and $D_2$, which we name so that the
orientation on c is compatible with that of $D_1$ and
incompatible with that of $D_2$.  We define F on each disk as
follows:

\begin{equation}
F(x) = \begin{cases} \pi_1 \circ F_0& \text{if x is in
$D_1$,}\\
\pi_2 \circ F_0&\text{if x is in $D_2$.}
\end{cases}
\end{equation}

F is continuous because $F_0(c)$ lies outside of U, and $\pi_1 =
\pi_2$ outside U.  F is degree 1 because $F_0: c \rightarrow \mathbb{R}^2
- T$ is degree 1.  F has Lipshitz constant $(2 + \sqrt2) \pi /
D$.

If the connected component of $S(p,R)$ containing q and r is not
path connected, it is still true that its $\delta$-neighborhood
is path connected for every $\delta > 0$, and we can take
$\gamma_{qr}$ to be a path lying as close as we like to $S(p,R)$. 
This curve $\gamma_{qr}$ suffices for our construction, producing
a map with Lipshitz constant $\epsilon + (2 + \sqrt2) \pi /D$ for
an arbitrary positive $\epsilon$.  Since the space of maps with a
fixed Lipshitz constant between compact spaces is compact, there
is a limit of these maps which is Lipshitz with constant $(2 +
\sqrt2) \pi / D$.

To prove the lower bound, we recall a theorem of Gromov about the
Uryson width.  In section F of \cite{G1}, Gromov proves that if a
Riemannian manifold $(M,g)$ has Uryson 1-width less than $\pi/2$,
then any contracting map from $(M,g)$ to the unit 2-sphere is
null-homotopic.  (In section 5, we discuss this result of Gromov
in detail.)  The quotient map from $(S^2, g)$ to the space of
connected components of the distance spheres around p has fibers
of diameter no more than D.  If this space of connected
components were a 1-polyhedron, then the Uryson 1-Width of
$(S^2,g)$ would be no more than D.  Applying Gromov's theorem, we
would see that the Lipshitz constant of a degree 1 map from
$(S^2, g)$ to the unit 2-sphere is at least $\pi / (2D)$. 

We deal with the general case by perturbing $dist(\cdot, p)$ to a
Morse function.  The space of connected components of the fibers
of a Morse function is an honest 1-polyhedron.  It is an exercise
in point-set topology to check that for each $\epsilon$ there is
a $\delta$ so that any function g obeying $|g(x) - dist(x,p)| <
\delta$ has fibers whose connected components have diameter less
than $D+ \epsilon$.  We take a Morse function g obeying this
inequality, which proves that the Uryson 1-Width of $(S^2, g) <
D+\epsilon$ for every $\epsilon$. \endproof

As a corollary, we are able to estimate the degree 1
hypersphericity, the hypersphericity, and the Uryson 1-Width of
$(S^2, g)$ in terms of the number D.  

$$D / (\pi (2 + \sqrt2)) \le \hbox{degree 1 HS}(S^2, g) \le
\hbox{HS}(S^2, g) \le 2 \hbox{UW}_1(S^2,g) / \pi \le 2D / \pi.$$

In this formula, $\hbox{degree 1 HS}$ stands for the degree 1
hypersphericity, $\hbox{HS}$ stands for the hypersphericity, and
$\hbox{UW}_1$ stands for the Uryson 1-Width.  In particular, we
see that the degree 1 hypersphericity, the hypersphericity, and
the Uryson 1-Width of a 2-sphere all agree up to a bounded
factor.  By contrast, we will prove in section 5 that for
surfaces of high genus, any two of these invariants can disagree
by an arbitrarily large factor.

\section{Families of Short Curves and Closed Geodesics}

In this section we will bound the length of the shortest closed
geodesic on $(S^2, g)$ by a multiple of the degree 1
hypersphericity.  The existence of a short closed geodesic
follows from fairly standard methods once we have a sufficiently
nice family of short curves on $(S^2, g)$.  We begin by proving
the existence of such a family.

\begin{theorem} Suppose $(S^2, g)$ has degree 1 hypersphericity
less than 1.  Then there is a map F from $(S^2, g)$ to a
trivalent tree T, so that each fiber of the map has length less
than 120.  The fibers of the map have controlled topology: the
preimage of any point in an edge of T is a circle, the preimage
of a terminal vertex of T is a point, and the preimage of a
trivalent vertex of T is homeomorphic to the greek letter
$\theta$.
\end{theorem}

Here is an outline of the proof.  The first step is to find a
collection of short disjoint circles in $(S^2, g)$, cutting it
into components of bounded diameter.  We begin by taking the
connected components of distance spheres around some point p in
$(S^2,g)$.  These disjoint circles have bounded diameter, but
they may be very long.  We replace each distance sphere by a
homologous curve which is a union of short circles.  These curves
are short, but they may not be disjoint.  By applying an
appropriate minimization procedure to this set of circles, we make
them disjoint.  These short disjoint circles will be fibers of
the map that we are going to construct.  The second step of the
proof is to fill in the map on each region between the short
disjoint circles.  Since these regions have small diameter, all
the geodesic segments in them are short.  We construct a map on
each region whose fibers consist of two or three geodesic
segments plus some curves of negligible length.  As a result, all
of the fibers have bounded length.

\proof Pick a point p in $(S^2, g)$, and consider
the distance spheres around p.  (We first smooth the distance
function so that almost every distance sphere is a union of
circles.)  Let $S_i$ be the distance sphere around p of radius
$12 i$, and let $T_i$ be the distance sphere around p of radius
$12i + 6$.  Of course $T_i$ separates $S_i$ from $S_{i+1}$. 
Equivalently, each curve from $S_i$ to $S_{i+1}$ has non-zero
(topological) intersection number with one of the components of
$T_i$.  By Theorem 1.1, each component $T_{i, j}$ of $T_i$ has
diameter less than 12.

We now use a standard trick to replace each $T_{i, j}$ with a
union of shorter curves, which still separate $S_i$ from
$S_{i+1}$.  For each component $T_{i, j}$, take a sequence of
points $q_{i, j, k}$ going around it, with consecutive points
very close together.  Let $\gamma_{i, j, k}$ be a minimal
geodesic from $q_{i, j, k}$ to $q_{i, j, 0}$.  Note that each
$\gamma_{i, j, k}$ is disjoint from $S_i$ and from $S_{i+1}$. 
Let $T_{i, j, k}$ be the closed curve formed by composing
$\gamma_{i, j, k}$, $\gamma_{i, j, k+1}$, and the portion of
$T_{i, j}$ from $q_{i, j, k}$ to $q_{i, j, k+1}$.  Now each
$T_{i, j, k}$ has length less than 24, each $T_{i, j, k}$ is
disjoint from $S_i$ and $S_{i+1}$, and each curve joining $S_i$
to $S_{i+1}$ has non-zero topological intersection number with
one of the $T_{i, j, k}$.  Unfortunately, the $T_{i,j,k}$ may be
far from disjoint.  One of the main steps of this proof is to
modify them to make them disjoint.

We will modify the $T_{i, j, k}$ using a minimization procedure. 
First, we define a minimal basis.  Let $(U,g)$ be a compact
surface with a Riemannian metric and convex boundary, and suppose
$H_1(U, \mathbb{Q})$ has dimension n.  We consider all sets of n
circles which form a rational basis for the first homology group
of U.  We define a partial ordering on these bases as follows. 
Let S be a set of n oriented curves in U, which form a rational
basis for $H_1(U)$, and let the lengths of the curves be $s_1 \le
... \le s_n$.  For two bases S and T, we say that T is smaller
than S if $t_i \le s_i$ for every i and $t_i < s_i$ for some i. 
We call S a minimal basis if there is no smaller basis.

Because U is compact, there is a minimal basis.  Each curve in a
minimal basis minimizes length in its homology class.  Because
the boundary of U is convex, each curve in a minimal basis is a
closed geodesic.  Also, no curve is a multiple covering of a
closed geodesic, because we could replace it with a single
covering of that geodesic, making the basis smaller.  These
geodesics must be distinct, else they would not form a basis, and
so they intersect transversally.

Following Gromov, we define a curve $\gamma$ to be straight if
for any two points x and y on $\gamma$, the distance in U from x
to y is equal to the distance in $\gamma$ from x to y.  In
section 5 of \cite{G2}, Gromov showed that a minimal rational
basis for the first homology group consists of straight circles. 
We will prove a slightly stronger statement.  We define a curve
$\gamma$ to be strictly straight if, for any two points x and y on
$\gamma$ and any minimal segment p from x to y, p lies in
$\gamma$.  A strictly straight circle is always an embedded
closed geodesic.

\begin{lemma} If $(U,g)$ is a compact surface with convex
boundary and S is a minimal rational basis for the first homology
of U, then every curve in S is strictly straight.  Any two curves
in S intersect at most once.  If U is planar, the curves in S are
disjoint.
\end{lemma}

\proof Suppose that one curve $\gamma$ in S is not strictly
straight.  Then there are points x and y in $\gamma$ and a
minimal segment p from x to y which does not lie in $\gamma$. 
Since p and $\gamma$ are both geodesics, p must be transverse to
$\gamma$ at x and y.  Now, we define $\gamma_1$ and $\gamma_2$ to
be the closed curves formed by connecting p to either path in
$\gamma$ from x to y.  Since p is minimal, each of $\gamma_1$ and
$\gamma_2$ is no longer than $\gamma$.  Moreover, $\gamma_1$ and
$\gamma_2$ have corners at x and y, so they can be homotoped to
curves which are strictly shorter than $\gamma$.  Now either
$\gamma_1$ or $\gamma_2$ is linearly independent from the other
curves in the basis S.  Replacing $\gamma$ in S by one of these
curves gives a strictly smaller rational basis.  Therefore, every
curve in a minimal basis is strictly straight.

Suppose that two geodesics $\gamma_1$ and $\gamma_2$ in S
intersect in two points x and y.  Let p be any minimal segment
from x to y.  Since $\gamma_1$ and $\gamma_2$ are strictly
straight, p lies in both $\gamma_1$ and $\gamma_2$, and so
$\gamma_1$ and $\gamma_2$ coincide.  Since $\gamma_1$ and
$\gamma_2$ are distinct, they do not intersect in two points. 

Now suppose that U is a planar domain.  Suppose that two of the
curves in S intersect in only one point.  Since they are
distinct closed geodesics, they intersect transversally, and
their intersection number must be 1 or -1.  But since U is a
planar domain, the intersection number of any two closed curves
in U is zero.  Therefore, the curves in S are disjoint. \endproof

With this lemma, we can modify the curves $T_{i,j,k}$ to make
them disjoint.  Let $U_i$ be the portion of $(S^2, g)$ between
$S_i$ and $S_{i+1}$.  By making an arbitrarily small bilipshitz
change of the metric, we can assume that each $U_i$ has convex
boundary.  To see this, for instance, we can use the normal
exponential map to give a metric with a geodesic boundary and
then slightly increase the metric along the boundary.  It is not
difficult. Then, in each $U_i$, we consider a minimal rational
basis of the first homology.  Notice that each $U_i$ is a planar
domain.  By Lemma 2.2, a minimal basis is realized by disjoint
curves.  We keep only the curves in the basis with length less
than 24, and label these $\gamma_{i,j}$.  Since each $T_{i,j,k}$
has length less than 24, the homology class of each $T_{i,j,k}$
lies in the rational span of the $\gamma_{i,j}$.  Now any path
from $S_i$ to $S_{i+1}$ has non-zero intersection number with
some $T_{i,j,k}$, and so it must have non-zero intersection
number with some $\gamma_{i,j}$.  Therefore, the $\gamma_{i,j}$
separate $S_i$ from $S_{i+1}$.

Cut $(S^2, g)$ along the curves $\gamma_{i,j}$.  We claim that
each component of the complement has radius less than 36.  Let X
be a component of the complement, not containing p, and B the
boundary component of X closest to p.  Any point in X is within
24 from B, using the distance function of $(S^2, g)$.  This
estimate remains true for the intrinsic distance in X.  To see
this, pick any point x in X, and consider a minimal geodesic g
from x to p.  We break g into pieces $g_i = g \cap U_i$.  Each
curve $g_i$ is a single minimal geodesic in $U_i$.  Then we
replace each $g_i$ with the minimal geodesic $\tilde g_i$ in the
slightly modified metric on $U_i$ between the endpoints of $g_i$. 
We define $\tilde g$ to be the union of all $\tilde g_i$.  The
curve $\tilde g$ still goes from x to p, and it has practically
the same length as g.  Since each $\gamma_{i,j}$ is strictly
straight in the modified metric on $U_i$, $\tilde g_i$ crosses
$\gamma_{i,j}$ no more than once.  Therefore, $\tilde g$ crosses
each $\gamma_{i,j}$ no more than once.  In particular, $\tilde g$
stays in X until it hits B, and we see that the intrinsic
distance from x to B is less than 24.  Now pick any point b in B. 
The distance from any point in X to b is less than 36.  On the
other hand, if X contains p, and x is any point in X, then the
curve $\tilde g$ from x to p stays in X and has length less than
24.

Let us summarize what we have done so far.  We have found a set
of disjoint embedded closed curves $\gamma_i$ on $(S^2, g)$, each
of length less than 24, so that each component of their
complement has radius less than 36.  The curves $\gamma_i$ which
we have just constructed will be fibers of the map we are going
to build.  It remains for us to define the map on each component
of their complement.  We call the components $U_k$.  For each
$U_k$, we will construct a map $F_k$ from $U_k$ to a tree $T_k$,
with the property that each boundary component of $U_k$ is the
fiber of a terminal vertex of $T_k$.  Putting together all of
these maps, we get a single map from $(S^2, g)$ to a tree T whose
fibers are just the fibers of the maps $F_k$.  Therefore, our
theorem reduces to the following lemma.

\begin{lemma} Let U be a planar domain, compact with boundary. 
Suppose each boundary component of U has length less than 24, and
that the radius of U is less than 36.  Then U admits a map to a
trivalent tree T, so that each fiber has length less than 120. 
Each boundary component of U is the fiber of a terminal vertex of
T.  The fibers of the map have controlled topology: the preimage
of a point on an edge is a circle, the preimage of a terminal
vertex is either a point or a boundary component, and the
preimage of a trivalent vertex is homeomorphic to the letter
$\theta$.
\end{lemma}

\proof We pick a point p so that every point in U lies in the
ball around p of radius 36.  After a small bilipshitz change of
the metric, we can assume that the boundary is convex. 
Therefore, every minimal path from an interior point to p avoids
the boundary and is a geodesic.  The proof will focus on the set
of all minimal geodesics in U with one endpoint at p.  Away from
the point p and from the cut locus of the exponential map, this
set of minimal geodesics forms a foliation, and we will view the
set of minimal geodesics on the whole space as a foliation with
some singular leaves.  The plan of the proof is to modify this
foliation until it consists of disjoint closed leaves (with mild
singularities) which will be the fibers of the desired map.

For an arbitrary metric, the foliation by minimal geodesics may
have very complicated singularities.  Morally, we should be able
to reduce the singularities to standard types by making a
$C^\infty$-small perturbation of the metric.  For technical
reasons, it is more convenient to perturb the geodesics in a
different way.  Let S be the circle around p of radius
$\epsilon$, where $\epsilon$ is much less than the injectivity
radius of the exponential map at p.  The normal exponential map
is a mapping from the normal bundle of S to U, and the images of
the fibers of the normal bundle are exactly the geodesics through
p.  Instead of perturbing the metric of U, we let $\tilde S$ be a
$C^\infty$-small perturbation of the circle S and use the normal
exponential map of $\tilde S$ in place of the normal exponential
map of S.

We now define our perturbed foliation precisely.  We parametrize
$\tilde S$ by the coordinate $\theta$.  For $r > 0$, we define
the normal exponential map $F(\theta, r)$ to be the endpoint of
the geodesic segment of length $r$ leaving $\tilde S$ at $\theta$
in the outward normal direction.  Our perturbed geodesics are
given by $F(\theta, r)$ for fixed $\theta$ and varying $r$.  The
minimal geodesics of our perturbed family extend only until the
infimal value of $r$ at which $F(\theta, r) = F(\theta', r')$
with $r' < r$.  This defines a foliation with singular leaves on
the exterior of $\tilde S$.

For each value of r, the function $F(\theta, r)$ defines a
Legendrian curve in U, and so we can view F as a 1-parameter
family of Legendrian curves.  The singularities of the foliation
by geodesics normal to $\tilde S$ correspond exactly to the
singularities of this family of Legendrian curves.  Because
$\tilde S$ was generic, this family contains only standard types
of singularities, as described in \cite{A}.  A generic Legendrian
curve has only two kinds of singularities: self-intersections and
cusps.  A generic 1-parameter family of Legendrian curves has
several more kinds of singularities, which occur at discrete
values of the parameter.  These singularities are
self-tangencies, triple-intersections, cusp formations and
cancellations, and cusp crossings. 

The possible singularities of the foliation are reduced further
by some metric considerations.  The Legendrian curve $F(\theta,
r)$ includes the set of all points at distance r from $\tilde S$,
and it may also include some points which are closer to $\tilde
S$.  We call the set of points at distance r from $\tilde S$ the
metric front of the Legendrian curve.  Because we are taking the
foliation by minimal geodesics, only the singularities on the
metric front of the Legendrian curve correspond to singularities
of the foliation.  Cusps, cusp cancellations, and cusp crossings
cannot occur on the metric front of the Legendrian curves
$F(\theta, r)$.  Therefore, the foliation by minimal geodesics
normal to $\tilde S$ only has singularities corresponding to
self-intersections, self-tangencies, cusp formations, and
triple-intersections.

A self-intersection of the Legendrian curve corresponds to a
double point of our foliation where two equal rays meet and
terminate.  Because self-intersections occur generically on
Legendrian curves, the double points form a 1-dimensional
manifold.  All other singularities occur at isolated points.  A
self-tangency of the Legendrian curve corresponds to a special
double point, where two equal rays meet head-on (making an angle
of $\pi$).  We will not need to distinguish these self-tangency
points from ordinary double points.  A cusp formation of the
Legendrian curve corresponds to a singular point of the foliation
where a single ray ends at a point.  A triple-intersection of the
Legendrian curve corresponds to a triple point of the foliation
where three equal rays meet and terminate. 

For now, we extend the foliation to the interior of the curve
$\tilde S$ by taking minimal geodesics to p.  The resulting
foliation on the interior of $\tilde S$ is non-singular except at
p, where it has a serious singularity which we will have to deal
with later.  By putting the boundary in general position, we can
assume that the only singular points which occur on the boundary
are double points.  There are a finite number of boundary double
points.  Since the boundary is convex, the foliation is
transverse to the boundary away from the boundary double points.

To summarize, we have constructed a foliation of U by perturbed
geodesics.  Because the perturbation was very slight, we may
assume that each perturbed geodesic has length less than 36.  The
foliation has five kinds of singular points: the point p, double
points, cusp formation points, triple points, and boundary double
points.  We include a diagram illustrating the five types of
singular points.

\includegraphics{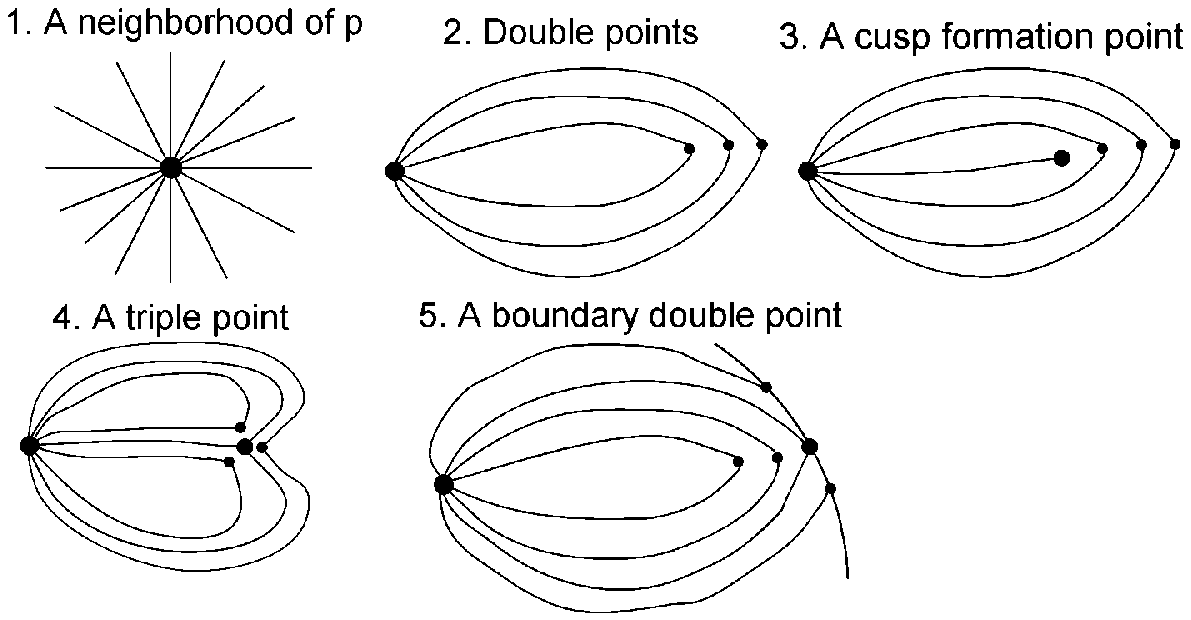}

We will do three surgeries to the foliation.  The first surgery
fixes the foliation so that each boundary component of U is a
leaf of the foliation.  It takes place in a small neighborhood of
the boundary of U.  The second surgery fixes the foliation so
that it extends over p.  It takes place in a small neighborhood
of p.  The second surgery uses the fact that U is a planar
domain.  The third surgery, which is only cosmetic, replaces each
singular leaf ending in a cusp formation point by a singular leaf
consisting of a single point.

We describe the first surgery.  Let c be a boundary component
of U.  The goal of the first surgery is to modify the foliation
near c so that c becomes a leaf.

Step 1. Modify the foliation in a neighborhood of each boundary
double point so that it is non-singular and transverse to c. 
(This creates a new triple point of the foliation in the interior
of U near the old boundary double point.) 

\includegraphics{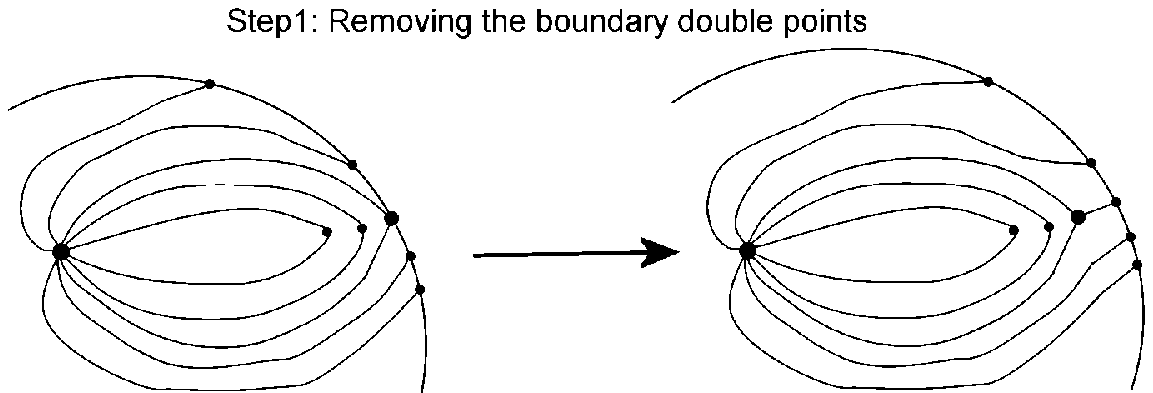}

Step 2. Parametrize a small neighborhood N of c by $S^1
\times [- \epsilon, \epsilon]$, so that each curve $\{ \theta \}
\times [-\epsilon, \epsilon]$ is a leaf of the foliation and each
circle $S^1 \times \{ r \}$ has length less than 24.  (The circle
with r-coordinate $- \epsilon$ is the boundary c.)  Let all the
circles with r-coordinate less than 0 be leaves of the modified
foliation.  In particular, c is a leaf of the modified foliation. 
Next, make a singular leaf of the modified foliation by joining
the circle with r-coordinate 0 to two rays from that circle to p. 
The two rays should be chosen so that they bisect the circle. 

This singular leaf cuts the remainder of N into two components,
which are each parametrized by $(0, \pi) \times (0, \epsilon]$. 
On each component we will replace the radial foliation whose
leaves are $\{ \theta \} \times (0, \epsilon)$ by a rectangular
foliation.

For each number $a$ strictly between 0 and 1/2, we make a leaf of
the rectangular foliation consisting of two radial line segments
and one circular arc.  The radial line segments have
$\theta$-coordinate equal to $\pi a$ and $\pi (1 - a)$ and
r-coordinate greater than or equal to $a \epsilon$.  The circular
arc has r-coordinate equal to $a \epsilon$ and
$\theta$-coordinate ranging from $\pi a$ to $\pi (1 - a)$. 
Finally, there is an exceptional leaf consisting of a radial line
with $\theta$-coordinate $\pi/2$ and r-coordinate greater than or
equal to $\epsilon/2$.  The foliation has a singularity at
$(\pi/2, \epsilon/2)$, where it is homeomorphic to a cusp
formation point.  This finishes the first surgery.

\includegraphics{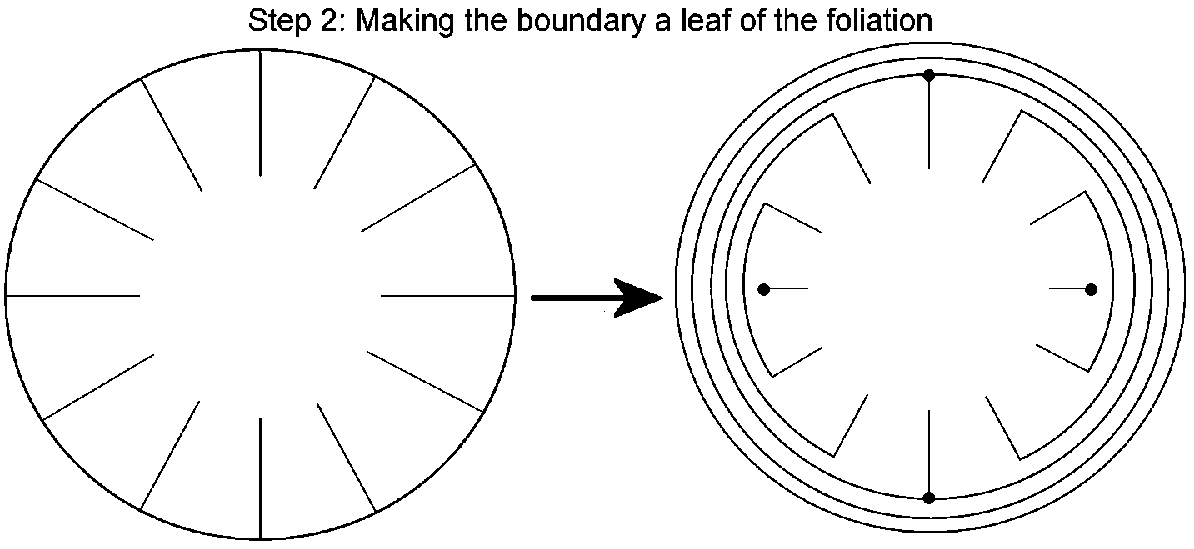}

At the end of the first surgery, all but finitely many leaves of
the foliation are curves that start and end at p and have length
less than 84.  There are three kinds of singular leaves.  First,
there are leaves with one triple point, which we will call triple
leaves.  The first surgery turned each boundary double point into
a triple point.  Each new triple point lies on a leaf consisting
of three rays and an arc parallel to the boundary component c. 
Since this arc has length less than 12, the length of each new
triple leaf is less than 120.  In addition, the original triple
points are still present.  Each original triple point lies on a
leaf consisting of three rays, with total length less than 108. 
Second, there are leaves with two rays joining a circle which is
close to a boundary circle, which we will call boundary leaves. 
Since each boundary circle has length less than 24, each boundary
leaf has length less than 96.  Finally, there are rays ending at
cusp formation points, which we will call cusp leaves.  Every
leaf has length less than 120.

We describe the second surgery.  The goal of the second surgery
is to modify the foliation near p so that it extends over p.

Step 1. Draw a small circle c around p, and foliate
the disk bounded by c by small concentric circles, with a
point at p.  Pick two points on c whose rays meet at a double
point, and add the rays to the circle to make a singular leaf
L.

\includegraphics{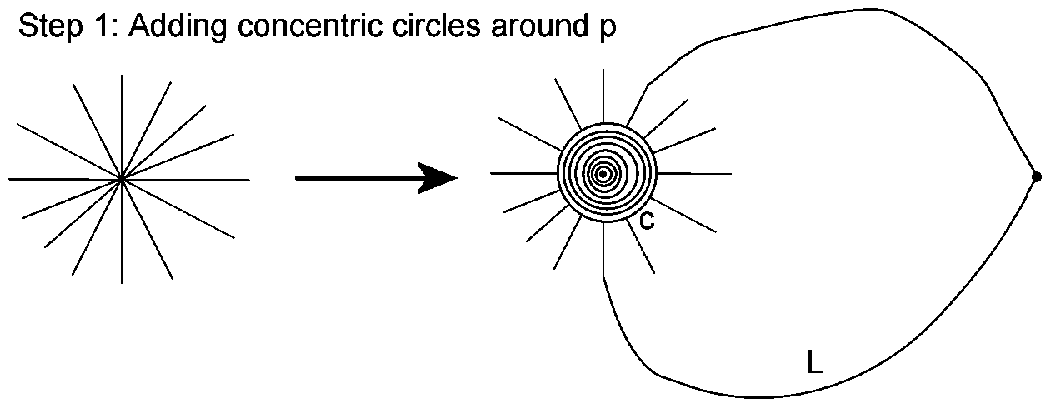}

Because U is a planar domain, the leaf L divides the remainder of
U into two components, which we will repair one at a time.  Each
component is a planar domain bounded by an arc a of the circle c
and a double leaf L joining the two endpoints of a.  We have a
foliation on the interior of this planar domain which extends
smoothly to L but which is transverse to a, and we need to modify
it to make it tangent to a.

Step 2. We define S to be the set of all points where a singular
leaf meets the arc a.  If the set S were empty, then each point
of a would lie on a double leaf.  The map taking each point to
the point on the other side of its double leaf would be a
continuous involution of the arc a with no fixed points, which is
impossible.  If there is exactly one point in S, then V is a disk
containing one cusp leaf and no other singular leaves.  In this
case, we are in a situation similar to the last part of the first
surgery.  We parametrize $a$ by $(0,1)$, so that the cusp leaf
crosses $a$ at 1/2, and so that for each number x between 0 and
1, the leaf crossing $a$ at x meets the leaf crossing $a$ at
$1-x$ in a double point.  We then replace the foliation in a
neighborhood of $a$ by a rectangular foliation, as in the first
surgery.

\vskip5pt
\includegraphics{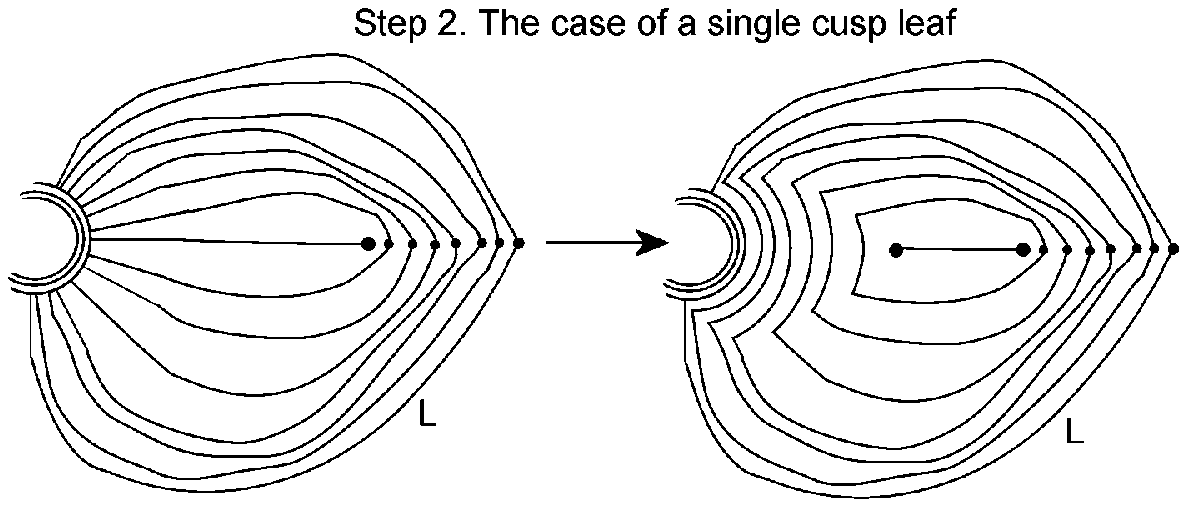}

If S has more than one point, then we consider the two outermost
points of S.  By a continuity argument, these two points must
belong to the same leaf, which we call L'.

The foliation between L and L' does not contain any singular
leaves and so it must have a very simple form.  Between L and L',
there are two intervals of a.  We parametrize the first interval
by $(0, 1/2)$ and the second interval by $(3/2, 2)$ so that the
leaf crossing $a$ at x meets the leaf crossing $a$ at 2-x in a
double point.  We let $a'$ be a curve very close to a but
slightly on the interior of the planar domain (slightly farther
from p).  We then replace the foliation between $a$ and $a'$ by a
rectangular foliation as in the first surgery.

\includegraphics{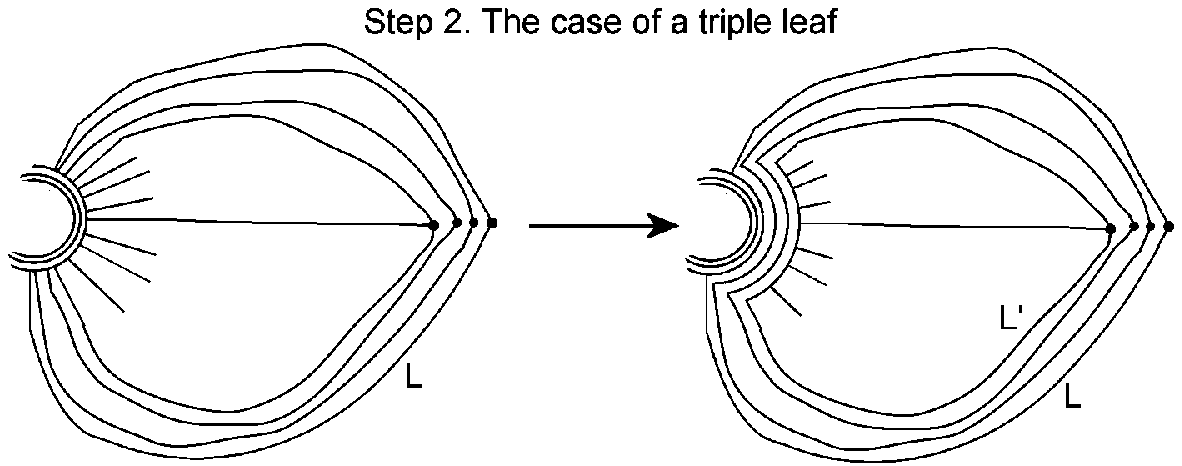}

We have now extended our foliation a bit further.  We have not
yet repaired the region between $L'$ and $a'$.  If $L'$ is a
triple leaf, this region has two components, and if $L'$ is a
boundary leaf, this region has one component.  Each component is
a planar domain bounded by an arc a' and a double leaf L' joining
the two endpoints of a'.  We have a foliation on the interior of
this planar domain which extends smoothly to L' but which is
transverse to a', and we need to modify it to make it tangent to
a'.  Since these are exactly the assumptions we made at the
beginning of Step 2, we can continue our construction
inductively, finding the outermost singular leaf of each new
component, and extending the foliation up to that leaf.  Since
there are only finitely many singular leaves, this process
terminates.

At the end of the surgery, the foliation is no longer singular at
the point p.  Since the circles and arcs in the surgery can be
taken arbitrarily small, the surgery has a negligible effect on
the lengths of leaves.  There are now only two kinds of singular
leaves.  The triple leaves and boundary leaves are both
homeomorphic to the letter $\theta$, and the cusp leaves are
homeomorphic to line segments with one cusp formation singularity
at each end.

We describe the third surgery.  The goal of the third surgery is
to modify the foliation in a small neighborhood of a cusp leaf to
replace the cusp leaf by a point. 

In the neighborhood of a cusp leaf, the foliation is homeomorphic
to the level sets of the distance function from a line segment in
$\mathbb{R}^2$.  The distance function admits a compactly
supported modification so that it has only one non-degenerate
minimum and no other criticial points.  The level sets of the
modified function give a new foliation with a singular point
instead of a singular line segment.  This operation increases the
lengths of leaves by an arbitrarily small amount.

Finally, we define the target tree to be the space of leaves of
our foliation and the map to be the quotient map.  All the
claims of the lemma have been satisfied. \endproof

This finishes the proof of Theorem 2.1.  Using this theorem, we
can bound the length of a closed geodesic.

\begin{theorem} Suppose $(S^2, g)$ has degree 1 hypersphericity
less than 1.  Then it contains a closed geodesic of length less
than 160.
\end{theorem}

\proof We will use two well-known existence theorems for closed
geodesics, which go back to Birkhoff.  The first result is the
curve shortening theorem.  Starting with any closed curve c in
$(S^2, g)$ which is not a closed geodesic, there is a canonical
homotopy which decreases the length of the curve over time and
which either converges to a closed geodesic or contracts the
curve to a point.  The second result concerns families of short
curves.  If there is a family of curves in $(S^2, g)$ of length
less than L which sweep out a degree non-zero map from $S^2$ to
$(S^2, g)$, then there is a closed geodesic in $(S^2, g)$ of
length less than L.  This theorem can be proved by using Morse
theory on the space of maps from the circle to $(S^2, g)$.  Both
results are explained in Croke's paper \cite{C}.

We will use the family of short curves guaranteed by Theorem 2.1. 
In the special case that our tree has no triple points, then the
fibers of our map sweep out $(S^2, g)$ by circles of length less
than 120, and by Birkhoff's theorem there is a closed geodesic of
length less than 120.  The main point of this proof is to deal
effectively with the singular fibers.

It turns out that singular fibers homeomorphic to the wedge of
two circles are more convenient for our proof than singular
fibers homeomorphic to the letter $\theta$.  By modifying our map
near the singular fibers, we can replace the singular fibers by fibers
homeomorphic to the wedge of two circles.  

\includegraphics{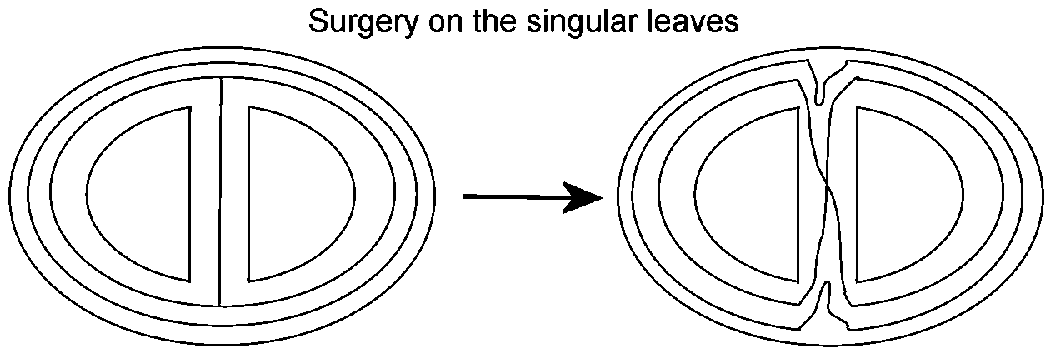}

After this modification, we have a map from $(S^2, g)$ to the
same tree, the preimage of each trivalent vertex is a wedge of
two circles, the preimage of each point on an edge of the tree is
a circle, and the preimage of each terminal vertex is a point. 
Each fiber has length less than 157.

We will use a curve shortening lemma for wedges of two circles,
which we now prove.  Let f be a map from the wedge of two circles
to $(S^2, g)$ with length less than 157.  We will prove that f
can be continuously deformed within the space of Lipshitz maps
until it reaches a local minimum of the length functional, and so
that every map on the way has length less than 160.  We can
choose a constant $\epsilon$ sufficiently small that if $f_0$ and
$f_1$ are two maps from the wedge of two circles to $(S^2, g)$,
and if the $C^0$ distance between them is less than $\epsilon$,
then there is a homotopy of maps $f_t$ between them so that the
length of each $f_t$ is less than 160.  Consider all the
1-parameter families of maps starting with f and with length less
than 160.  Let l be the infimal length of all maps in such
families.  Let $f_i$ be a sequence of maps in such families with
the length of $f_i$ equal to $l_i \rightarrow l$.  We can assume
that each $f_i$ has length less than $157$.  After
reparametrization, we can assume that each $f_i$ has Lipshitz
constant less than 1000.  The space of maps with Lipshitz
constant less than 1000 is precompact in the $C^0$ topology, so
the $f_i$ have a $C^0$ sublimit $f_{\infty}$, and the length of
$f_{\infty}$ is less than or equal to l.  For large i, the $C^0$
distance from $f_i$ to $f_{\infty}$ is less than $\epsilon$, and
so we can deform $f_i$ to $f_{\infty}$ without increasing its
length over 160.  Thus we can deform f to $f_{\infty}$ without
increasing its length over 160, and the length of $f_{\infty}$
must be equal to l.  If any function g was within $\epsilon$ of
$f_{\infty}$ in the $C^0$ metric and had length less than l, we
could deform $f_{\infty}$ to g without increasing its length over
160, contradicting the definition of l.  Therefore, $f_{\infty}$
minimizes length in a $C^0$ neighborhood.

We now describe the possible maps of a wedge of two circles into
a surface which are local minima for the length function.  We let
p be the image under a map f of the base point of the wedge of
two circles.  Clearly, each circle in the wedge is mapped to a
geodesic starting and ending at p, which might or might not be
the constant map.  If exactly one of these circles is constant,
then the other one must be a closed geodesic.  If neither circle
is constant, then there are four geodesic segments leaving p, and
the length-minimizing condition implies that the sum of the four
corresponding unit vectors is zero.  Because the tangent space at
p is 2-dimensional, this condition guarantees that the 4 vectors
consist of two pairs of opposite vectors.  Therefore, the image
of the wedge of two circles consists of either 1 or 2 closed
geodesics.  To summarize, either there is a closed geodesic of
length less than 160, or each singular fiber can be contracted to
a point through a family of wedges of two circles with length
less than 160.

From now on, we will assume that there is no closed geodesic of
length less than 160 and work up to a contradiction.  Because of
the standard curve-shortening theorem, we know that each
non-singular fiber contracts to a point through curves of length
less than 160.  Each family of curves contracting a non-singular
fiber to a point gives rise to a map from a disk to $(S^2, g)$,
the boundary of the disk mapping to the fiber.  We will refer to
a family of curves of length less than 160 which contracts a
non-singular fiber to a point as a contracting disk for that
fiber.  Any two contracting disks for the same fiber can be glued
together to get a map from the 2-sphere to $(S^2, g)$, and if the
degree of this map is non-zero then Birkhoff's theorem guarantees
a closed geodesic of length less than 160.  Therefore, for each
non-singular fiber, all the contracting disks are homologous. 

Pick one terminal vertex t of T as a root and let p be the
preimage of this vertex in $(S^2, g)$.  Each non-singular fiber
divides $(S^2, g)$ into two disks, and we will refer to the disk
that does not contain p as the far disk of that non-singular
fiber.  We will now prove that the contracting disk of every
non-singular fiber is homologous to its far disk.  We will work
inductively from the outermost branches of the tree towards the
root.  For fibers on the same edge as a terminal vertex other
than t, the fibering itself gives a contraction of the fiber to a
point whose contracting disk is the far disk of the fiber. 
Notice that once we have proved our claim for one fiber, then the
claim follows for all other fibers over the same edge of T. 
Therefore, it suffices to prove our claim inductively one edge at
a time.

We consider an edge $E$.  Let u be the endpoint of $E$ farthest from
t, and let $E_1$ and $E_2$ be the two edges meeting $E$ at u.  By
induction, we may assume that our claim holds for both $E_1$ and
$E_2$, and we will prove that it holds for $E$.  Let $F$ be the fiber
over a point in $E$ very close to u, and let $F_1$ and $F_2$ be
fibers over points in $E_1$ and $E_2$ very close to u.  Orient
the fibers so that $F$ is homologous to the sum of $F_1$ and $F_2$
within the inverse image of a small neighborhood of u.  We know
that the singular fiber over u contracts to a point through a
family of wedges of circles of length less than 160.  By
restricting in various ways, this one contraction produces a
contraction of $F$, $F_1$, and $F_2$.  We let $D$, $D_1$, and $D_2$
be the contracting disks for these fibers produced by the
contraction of the singular fiber.  Since these disks all come
from the same contraction of the singular fiber, $D$ is homologous
to $D_1 + D_2$.  By induction, we may assume that $D_1$ is
homologous to the far disk of $F_1$ and that $D_2$ is homologous
to the far disk of $F_2$.  Therefore $D$ is homologous to the
far disk of $F_1$ plus the far disk of $F_2$, which is the far
disk of $F$.

On the other hand, for fibers over the edge containing the root
t, the disk containing p is obviously a contracting disk.  Since
the far disk is also a contracting disk, Birkhoff's theorem
guarantees a closed geodesic of length less than 160. \endproof

\section{Applications to Surfaces in Euclidean Space}

In this section, we will prove three estimates for the
hypersphericity of a spherical or planar surface in Euclidean
space.  Recall that for a surface with boundary, $(\Sigma,
\partial \Sigma)$, the degree 1 hypersphericity is defined using
degree 1 maps from the pair $(\Sigma, \partial \Sigma)$ to the
pair $(S^2, *)$.  The first estimate that we will prove is a
slicing inequality.

\begin{theorem} (Slicing Inequality) Let f be an embedding of the
unit n-ball in $\mathbb{R}^n$, and let $f_0$ be the first n-2 coordinates
of f.  Then one of the level sets of $f_0$ contains a connected
component with degree 1 hypersphericity greater than 1/24.
\end{theorem}

The proof of this theorem is based on the following lemma.

\begin{lemma} Let $\pi$ be the projection from $S^n$ to the first
n-2 coordinates of $\mathbb{R}^{n+1}$.  Suppose that g is a Riemannian
metric on $S^n$, and that the degree 1 hypersphericity of every
fiber of $\pi$ is less than 1.  Then the Uryson (n-1)-width of
$(S^n, g)$ is less than 12.
\end{lemma}

\proof Let s be a continuous section of the projection $\pi$.  We
define the function $F(x)$ to be the distance, in the intrinsic
metric of the fiber through x, between x and the point
$s(\pi(x))$.  The function F is clearly continuous away from the
singular fibers of $\pi$.  Since the metric g is bounded by a
constant multiple of the standard metric on $S^n$, the function
F(x) goes to zero as x approaches any of the singular fibers. 
Therefore, F is continuous.  The function $(\pi, F)$ maps $S^n$ to
$\mathbb{R}^{n-1}$.  By Theorem 1.1, each connected component of
the fibers of the map $(\pi, F)$ has diameter less than 12.  We
consider the induced map from $(\pi, F)$ to the space of
connected components of level sets of $(\pi, F)$.  By definition,
the fibers of this map have diameter less than 12.  If the space
of connected components of level sets of $(\pi, F)$ were an
(n-1)-polyhedron, we would be done.

The space of connected components of level sets of $(\pi, F)$ may
be very complicated.  We explain how to approximate it by
(n-1)-dimensional polyhedra, by a method that goes back to
Alexandrov.  (For more information on polyhedral approximation,
see the book {\it Dimension Theory} by Hurewicz and Wallman,
\cite{HW}, especially section 9 of chapter 5.)  We cover
$\mathbb{R}^{n-1}$ by very small open sets
$B_i$.  We can arrange these open sets so that each point lies in
$B_i$ for at most n values of i.  The connected components of the
sets $(\pi, F)^{-1}(B_i)$ form an open cover of $(S^n,g)$.  There
may be infinitely many open sets in this cover, but since $S^n$
is compact we may take a finite subcover.  Each point of $(S^n,
g)$ lies in at most n of the open sets in the open cover.  We now
consider the map from $(S^n, g)$ to the nerve of this cover,
which is an honest (n-1)-dimensional polyhedron.  Each fiber of
this map lies in a connected component of the preimage of one of
the sets $B_i$.  By taking a very fine open cover, we can arrange
that the diameter of each connected component of the preimage of
each $B_i$ is less than 12. \endproof

Let us apply the above lemma to the pullback of the standard
metric on $S^n$ by an arbitrary diffeomorphism.  By Gromov's
estimate in \cite{G1}, the Uryson (n-1)-width of the unit
n-sphere is at least $\pi/2$.  Therefore, one of the level sets
of $\pi$ has hypersphericity at least $\pi / 24$.  Thus we
have proven a slicing inequality for the unit n-sphere. 

With these lemmas we give the proof of Theorem 3.1.  First,
rescale the function f so that the image lies in the unit ball of
$\mathbb{R}^n$.  This rescaling does not affect the fibers of
$f_0$.  Next, we embed the unit ball as a hemisphere of the unit
n-sphere, again without changing the first n-2 coordinate
functions.  Now, we put a metric g on the n-sphere as follows. 
On the image of the map f, we use the pushforward of the metric
on the unit ball.  On the complement of the image of f, we make
the Riemannian metric very small.  We are now in the situation of
Lemma 3.2.  Since the unit ball with the boundary collapsed to a
point admits a degree 1 map to the unit sphere with Lipshitz
constant $\pi$, the Uryson (n-1)-width of the metric g is at least
$1/2$.  Therefore, one of the fibers of the map $\pi$ has degree
1 hypersphericity at least $1/24$.  But the fibers of the map
$\pi$ are practically the level sets of $f_0$ with the boundary
collapsed to a point.  Therefore, one of the connected components
of a level set of $f_0$ has degree 1 hypersphericity at least
$1/24$.  This finishes the proof of Theorem 3.1.

Our next estimate is a variation of the isoperimetric inequality
in $\mathbb{R}^3$.

\begin{theorem} (Isoperimetric Inequality) Let U be a bounded
open set in $\mathbb{R}^3$ with boundary diffeomorphic to a 2-sphere. 
Then the boundary of U contains disjoint open sets $S_i$, so that
the following inequality holds.

$$\rm{Vol}(U) < 10^6 \sum \rm{HS}(S_i)^3.$$

In this equation, $\rm{HS}(S_i)$ stands for the degree 1
hypersphericity of $S_i$.
\end{theorem}

To get a sense of the estimate involved in this theorem, consider
how it applies to a rectangular solid R, with side lengths $R_1
\le R_2 \le R_3$.  The classical isoperimetric inequality gives a
bound for the volume of R on the order of $R_2^{3/2} R_3^{3/2}$. 
To apply our isoperimetric inequality, note that the boundary of
R has Uryson 1-width on the order of $R_2$.  Therefore, no set
$S_i$ has degree 1 hypersphericity greater than a multiple of
$R_2$.  Also, the degree 1 hypersphericity of $S_i$ is bounded by
a multiple of the square root of the area of $S_i$.  Putting
these bounds together, we see that $\sum HS(S_i)^3$ is bounded by
a multiple of $R_2^2 R_3$.  On the other hand, taking the sets
$S_i$ as disks of radius $R_2$, we get $\sum HS(S_i)^3$ on the
order of $R_2^2 R_3$.  Therefore, our isoperimetric inequality
gives a bound for the volume of R on the order of $R_2^2 R_3$. 
The constant in our isoperimetric inequality is much worse than
that in the standard inequality, but the dependence on $R_2$ and
$R_3$ is better.  Our inequality gives the best
dependence on $R_2$ and $R_3$ that can be expected, because the
boundary of the rectangle R can be isotoped without stretching to
enclose a volume on the order of $R_2^2 R_3$.

Now we turn to the proof of Theorem 3.3.  Fix a point p in the
interior of U.  For each number R, let B(R) be the ball of radius
R around p, and let $S(R)$ be the intersection of the boundary of
U with B(R).  We let $\phi_R$ be a diffeomorphism of
$\mathbb{R}^3$ with the unit 3-sphere minus a point Q, with the
following properties.  The diffeomorphism $\phi_R$ maps B(R) to
almost the entire unit 3-sphere using the exponential map at the
conjugate point to Q.  The diffeomorphism $\phi_R$ maps the
complement of B(R) into a tiny ball in the 3-sphere around the
point Q.  The Lipshitz constant of $\phi_R$ is $\pi/R$.  Using
the maps $\phi_R$, we define a map F from $\partial U \times (0,
\infty)$ to the unit 3-sphere, by taking $F(x,R) = \phi_R(x)$. 
For very large values of R, $\phi_R$ maps the boundary of U into
a tiny neighborhood of the conjugate point of Q, and for very
small values of R, $\phi_R$ maps the boundary of U into a tiny
neighborhood of Q.  Therefore, we can easily complete F to give a
degree 1 map from $\partial U \times [0, \infty]$ to the unit
3-sphere, mapping $\partial U \times \infty$ to the conjugate
point of Q and $\partial U \times 0$ to Q.

Now we assign a metric g to $\partial U \times [0,
\infty]$, so that for each R, the restriction of g to the fiber $\partial
U$ is given by taking the original metric on S(R), rescaling it
by $\pi/R$, and making the metric on the rest of the boundary of
U very small.  Then we make the metric g very large in the
directions transverse to the fibers.  In this way, we arrange
that the map F has Lipshitz constant 1, which shows that the
metric g has hypersphericity 1 and hence Uryson 2-width at least
$\pi/2$.  By Lemma 3.2., one of the fibers of $\partial U \times
(0, \infty)$ has degree 1 hypersphericity at least $\pi/24$. 
But each fiber is practically a rescaling of $S(R)$ with the
boundary contracted to a point.  Therefore, for some value of R,
$S(R)$ has degree 1 hypersphericity at least $R/24$.

For each point p in U, call $B(p,R)$ good if S(R) has degree 1
hypersphericity at least R/24.  By the Vitali covering lemma,
there are disjoint good balls $B(p_i, R_i)$ with total volume at
least one eighth the volume of U.  Now an easy computation shows
that $\sum \rm{HS}(S(p_i, R_i))^3 > \rm{Vol}(U) / 10^6$.
\endproof

The theorem that we just proved did not impose any lower bound on
the hypersphericity of the boundary of U.  In order to prove such
a bound, we have to make a stronger assumption about the geometry
of U.  In particular, we will prove the following variation of the
isoperimetric inequality in $\mathbb{R}^3$.

\begin{theorem} (Isoperimetric Inequality) Let U be a bounded
open set in $\mathbb{R}^3$ with boundary diffeomorphic to a 2-sphere. 
Suppose that there is an area-contracting map of non-zero degree
from the pair $(U, \partial U)$ to the unit 3-sphere $(S^3, *)$. 
Then the degree 1 hypersphericity of the boundary of U is at
least 1/60.
\end{theorem}

\proof The main ingredient of the proof is Theorem 2.1. 
According to Theorem 2.1, if the boundary of U has degree 1
hypersphericity less than 1/60, then it admits a map $\pi$ to a
trivalent tree T whose fibers have length less than 2.  The
fibers also have controlled topology, and in particular the fiber
of each point on an open edge of T is a circle.

Each of these fibers bounds a disk with bounded area in
$\mathbb{R}^3$.  These small area disks are useful in controlling
the degree of area-contracting maps.  Unfortunately, the small
disks need not lie in U and so they can intersect in a complicated way. 
To circumvent this problem, we will construct a Riemannian metric
g on the 3-ball, whose restriction to the boundary is isometric
to $\partial U$.  The metric on the interior will be sufficiently
large that for any two points x and y on the boundary, the
distance in $(B^3, g)$ between x and y is equal to the distance
in $\partial U$ between x and y.  At the same time, the metric on
the interior will be sufficiently small that it can be cut by
disjoint disks of bounded area into components of bounded volume.

We first choose a very large but finite collection of the fibers
of the map $\pi$.  Each fiber in our collection will be a circle. 
By choosing sufficiently many fibers, we may assume that each
component of their complement has very small area and has at most
three boundary components.  We call the fibers in our collection
$\{ F_i \}$.  At each fiber $F_i$, we attach a 2-cell $D_i$ to
$\partial U$ with an attaching map taking the boundary of the
2-cell diffeomorphically to $F_i$.  We give the 2-cell $D_i$ the
Riemannian metric of a hemisphere of radius $R_i$, where $2 \pi
R_i$ is the length of the fiber $F_i$, in such a way that the
attaching map is an isometry.  Since the length of $F_i$ is less
than 2, the radius of $D_i$ is less than $1/ \pi$ and the area
of $D_i$ is less than $2 / \pi$. 

For each component $A_j$ of the complement of the curves $F_i$ in
$\partial U$, we consider the sphere $S_j$ consisting of the
union of $A_j$ and the disks filling each component of the
boundary of $A_j$.  At each 2-sphere $S_j$, we attach a 3-cell
$B_j$ to our 2-complex with an attaching map taking the boundary
of the 3-cell diffeomorphically to $S_j$.  We give the 3-cell
$B_j$ a Riemannian metric consisting of a cylinder $S_j \times
[0, 1]$ capped by a very small metric on a ball joined to the
end $S_j \times \{ 1 \}$.  The attaching map is the identity
isometry from $S_j \times \{ 0 \}$ to $S_j$.  The resulting
3-complex is homeomorphic to a 3-ball with boundary $\partial U$,
and we have equipped it with a Riemannian metric which we will
call g.

Because the area of $A_j$ is negligibly small and because the
boundary of $A_j$ consists of at most three fibers, the sphere
$S_j$ has area less than $6 / \pi$.  The volume of $B_j$ is
negligibly more than the area of $S_j$ times $1$, so the volume
of $B_j$ is less than $6 / \pi$.  Next we show that the diameter
of $S_j$ is less than $2$.  Let x and y be any two points in
$S_j$.  Suppose for now that x is in $D_1$ and y is in $D_2$. 
Because we have chosen the fibers $F_i$ extremely densely, there
is a point z in $A_j$ which is extremely close to both $F_1$ and
$F_2$.  The distance in $S_j$ from z to either x or y is less
than $1$, and so the distance from x to y is less than $2$.  If x
and y are both in the same disk, then the distance between them
is less than $1$.  If one or both of x and y lies in $A_j$, the
distance between x and y is still less than $2$ because any point
in $A_j$ can be moved a negligibly small distance into one of the
disks $D_i$ in $S_j$.

We check that the distance in $(B^3, g)$ between two points x and
y in $\partial U$ is equal to the distance between x and y in
$\partial U$.  Let p be any path in $(B^3, g)$ between x and y. 
For any 3-cell $B_j$, let $p_0$ be a connected component of $p
\cap B_j$, with endpoints $t_1$ and $t_2$ in $S_j$.  We replace
$p_0$ by a minimal geodesic in $S_j$ from $t_1$ to $t_2$, which
we call $p_0'$.  The length of $p_0'$ is no more than the
diameter of $S_j$ which is $2$.  If $p_0$ intersects $S_j
\times \{ 1 \} \subset B_j$, then it has length at least $2$, and
so it is at least as long as $p_0'$.  On the other
hand, if $p_0$ does not intersect $S_j \times \{ 1 \}
\subset B_j$, then $p_0$ lies in the product $S_j \times [0, 1]$
and is at least as long as $p_0'$.  Repeating this argument for
each component of p meeting each ball $B_j$, we replace p by a
path p' of equal or smaller length lying on the union of
$\partial U$ and the disks $D_i$.  For any 2-cell $D_i$, let
$p'_0$ be a connected component of $p' \cap D_i$ with endpoints
$u_1$ and $u_2$ in $F_i$.  We replace $p_0'$ by the minimal path
in $F_i$ from $u_1$ to $u_2$, which we call $p_0''$.  Since this
minimal path is a portion of the equator of the hemisphere $D_i$,
$p_0'$ is at least as long as $p_0''$.  Repeating this argument
for each component of p' meeting each disk $D_i$, we replace p'
by a path p'' of equal or smaller length lying in $\partial U$. 
This proves the estimate for the distance.

Since $\partial U$ is embedded in $\mathbb{R}^3$, we have an
isometric embedding map i from $\partial U$ to $\mathbb{R}^3$. 
We view the map i as a map from the boundary of $(B^3, g)$ to
$\mathbb{R}^3$.  If we consider $\partial U$ as a metric space
using the distance function of $(B^3, g)$, the map i has Lipshitz
constant 1.  We now use a simple extension result for Lipshitz
maps, which goes back to McShane in \cite{Mc}.  Any Lipshitz map
from a subspace of some metric space to the real numbers admits
an extension to the whole metric space with the same Lipshitz
constant.  We extend each coordinate functon of the embedding i
to all of $(B^3, g)$, giving a $\sqrt3$-Lipshitz map from $(B^3,
g)$ to $\mathbb{R}^3$.

Now let F be any area-contracting map from $(U, \partial U)$ to
the unit 3-sphere, taking the boundary of U to the base point of
$S^3$.  Our aim is to prove that F has degree 0.  We first extend
F to all of $\mathbb{R}^3$ by mapping the complement of U to the
basepoint of $S^3$.  We will consider $F \circ i$, which is a map
from $(B^3, g)$ to the unit 3-sphere, taking the boundary of
$B^3$ to the basepoint of $S^3$, and with the same degree as F. 
We now use our decomposition of $B^3$.  Each 3-cell $B_j$ in
$(B^3, g)$ has volume less than $6 / \pi$, and so the image of
$B_j$ has volume less than $18 \sqrt{3} / \pi < 12$.  The
boundary of $B_j$ consists of $A_j$ together with at most three
disks $D_i$.  Each disk $D_i$ has area less than $2 / \pi$, and
so the image of each disk $D_i$ has area less than $6 / \pi$. 
Since $F \circ i$ maps $A_j$ to the basepoint of $S^3$, the image
of each disk is a 2-cycle.  Each 2-cycle in the unit 3-sphere of
area $A < 4 \pi$ bounds a 3-chain of volume less than $\pi^2 (A /
(4 \pi))^{3/2}$, so the image of each disk $D_i$ bounds a 3-chain
of volume less than $\pi^2 ((6/\pi) / (4 \pi))^{3/2} < 1$.  Let
$C_j$ be the 3-cycle in $S^3$ which is the union of the image of
$B_j$ and the 3-chains bounded by each $D_i$ in the boundary of
$B_j$.  Each $C_j$ has volume less than 15, and since the volume
of the unit 3-sphere is $2 \pi^2$, each $C_j$ is null-homologous. 
But the image of the fundamental homology class of $(B^3,
\partial B^3)$ is homologous to the union of the cycles $C_j$
with appropriate orientations.  Therefore, the map F has degree
0. \endproof

\section{Estimates for Surfaces of Non-Zero Genus}

In this section, we prove some estimates on the geometry of
surfaces of non-zero genus with bounded hypersphericity.

We first apply the methods of section 1 to estimate the systoles
of closed orientable surfaces.  In this paper, the systole of a
Riemannian surface will be the shortest length of a homologically
non-trivial closed curve in the surface.  (Sometimes systole refers to
the shortest homotopically non-trivial closed curve.)  To give
some context, we recall the main estimates for the systoles of
surfaces.  Besicovich proved that the systole of a surface of
area A is bounded by $\sqrt{A} / 2$, and Gromov proved that for
surfaces of high genus G the systole is bounded by $c(G)
\sqrt{A}$ where c(G) falls off almost as fast as $1 /
\sqrt{G}$. (See \cite{G2} for Gromov's theorem and chapter 4 of
\cite{G3} for some history and many other results.)  We will
bound the systole of a closed orientable surface by its
degree 1 hypersphericity.

\begin{theorem} Suppose that $(\Sigma, g)$ is a closed oriented
surface of genus at least 1, with degree 1 hypersphericity less
than 1.  Then $(\Sigma, g)$ contains a homologically non-trivial
closed curve of length less than $4 \pi$.
\end{theorem}

\proof Let $\gamma$ be a shortest homologically non-trivial curve in
$\Sigma$, and let L be the length of $\gamma$.  If L is less than
$4 \pi$, then we are done, so we may assume that L is at least $4
\pi$.  Pick two points
on $\gamma$, p and q, with $dist(p,q) = L/4$, and consider the
function $F_0(x) = (dist(p,x), dist(q,x))$.  As Gromov proved in
section 5 of \cite{G2}, the distance between two points of
$\gamma$ is the length of the shorter of the two paths in
$\gamma$ joining the two points.  (This result of Gromov is easy,
and a stronger easy result is proved in Lemma 2.2 of this paper.) 
Therefore, we can exactly determine the restriction of $F_0$ to
$\gamma$.  The image of $F_0$ is the square $S_1$ with vertices
$(0, L/4), (L/4, 0), (L/2, L/4),$ and $(L/4, L/2)$.

We will prove that any curve homologous to $\gamma$ includes a
point within $2 \pi$ of q.  Let $\tau$ be any curve homologous to
$\gamma$, and suppose that the distance from $\tau$ to q is at
least $2 \pi$.  In that case, $F_0(\tau)$ lies above the line $y
= 2 \pi$ in $\mathbb{R}^2$, and so avoids the square $S_2$, with vertices
$(L/4 - \pi, \pi), (L/4, 0), (L/4 + \pi, \pi)$, and $(L/4, 2
\pi)$.  Since $S_2$ is a square of side-length $\sqrt2 \pi$, it
contains a disk D of radius $\pi / \sqrt2$.  Since we assume $L >
4 \pi$, $S_2$ is contained in $S_1$, and therefore $F_0(\gamma)$
also avoids the interior of $S_2$.  Let $\Sigma_i$ be the
connected components of the complement of $\gamma$ and $\tau$ in
$\Sigma$.  (After putting $\tau$ in general position, we can
assume there are only finitely many connected components.)  Let
$B_i$ be the boundary of $\Sigma_i$.  Each $B_i$ consists of a
subset of $\gamma$ and a subset of $\tau$, so the map $F_0$ takes
$B_i$ to $\mathbb{R}^2 - D$.  Moreover, the portion of $F_0(B_i)$
below the line $y = 2 \pi$ must be either the portion of
$F_0(\gamma)$ below that line or else empty, and so the image of
$B_i$ has winding number -1, 0, or 1 around D.  Because $\gamma$
is homologous to $\tau$, an appropriate sum of integer multiples
of the $B_i$ is equal to $\gamma$ - $\tau$ (as 1-chains).  Since
the image of $\gamma$ - $\tau$ has winding number 1 around D, at
least one of the $B_i$ has non-zero winding number around D.  Say
that $B_1$ has winding number equal to 1 or -1.

We now define a contracting degree 1 map $F$ from $\Sigma$ to the
unit sphere.  Using the exponential map, we construct a
homeomorphism $\pi_1$ from the disk D to the Northern hemisphere
of the unit sphere.  We define a map $\pi_2$ from the disk D to
the Southern hemisphere of the unit sphere by composing $\pi_1$
with a reflection through the plane of the equator.  We extend
$\pi_1$ and $\pi_2$ to all of $\mathbb{R}^2$ as follows.  For any
point x outside of D, let p(x) be the nearest point to x on the
boundary of D.  Now, define $\pi_1(x) = \pi_1(p(x))$ and
$\pi_2(x) = \pi_2(p(x))$.  The maps $\pi_1$ and $\pi_2$ agree on
the the boundary of D and hence on all the points outside of D. 
The maps $\pi_1$ and $\pi_2$ have Lipshitz constant $1/\sqrt2$.
 
Finally we define our degree 1 map.  Recall that $\Sigma_1$ is an
open set in $\Sigma$ with boundary $B_1$.  Let $\Sigma_2$ be the
complement of $\Sigma_1$.  We define F on each component as
follows:

\begin{equation}
F(x) = \begin{cases} \pi_1 \circ F_0& \text{if x is in
$\Sigma_1$,}\\
\pi_2 \circ F_0&\text{if x is in $\Sigma_2$.}
\end{cases}
\end{equation}

F is continuous because $F_0(B_1)$ lies outside of D, and $\pi_1
= \pi_2$ outside D.  F is degree 1 or -1 because $F_0: B_1
\rightarrow \mathbb{R}^2 - D$ is degree 1 or -1.  After composing
with a reflection, we may assume that F is degree 1.  F has
Lipshitz constant $1$.  Since we assumed that the degree 1
hypersphericity of $(\Sigma, g)$ is less than 1, there is no such
map, and we conclude that the distance from $\tau$ to q is less
than $2 \pi$.

Since this estimate holds for any curve $\tau$ homologous to
$\gamma$, there is no curve homologous to $\gamma$ in the open
manifold $\Sigma - B(q, 2 \pi)$.  We claim that $B(q, 2 \pi)$
contains a curve with non-trivial homology class in
$H_1(\Sigma)$.  To see this, consider the Mayer-Vietoris sequence
associated to the covering $\Sigma = U \cup V$, where U is $B(q,
2 \pi)$ and V is a very small neighborhood of the complement of
$B(q,2 \pi)$: $H_1(U) \oplus H_1(V) \rightarrow H_1(\Sigma)
\rightarrow H_0 (U \cap V)$.

We know that the homology class of $\gamma$ is not in the image
of $H_1(V)$ in $H_1(\Sigma)$.  If the homology class of $\gamma$
is in the image of $H_1(U)$, we are done.  If not, then it has a
non-zero image $\alpha$ in $H_0(U
\cap V) = H_0(\partial U)$.  If $\partial U$ has k components,
$c_1, ...,
c_k$, then $H_0 (\partial U)$ has one generator $[c_i]$ for each
component $c_i$, and $\alpha = \sum \langle \gamma, c_i \rangle
[c_i]$, where $\langle \gamma, c_i \rangle$ is the intersection
number of $\gamma$ and $c_i$ in $\Sigma$. If $\alpha$ is not
zero, then one of the intersection numbers is not zero, and so
one of the boundary components is homologically non-trivial in
$\Sigma$.  In either case, we have found a homologically
non-trivial curve in $B(q, 2 \pi)$.  We can now factor that
curve, by adding many geodesics to and from q, into a sequence of
curves with length bounded by a number as close as we like to $4
\pi$.  One of these curves must have a non-trivial homology
class. \endproof

We now turn to two more difficult estimates on the geometry of a
surface of small hypersphericity.  The first estimate bounds the
lengths of G homologically independent curves on $\Sigma$.

\begin{theorem} Suppose that $(\Sigma, g)$ is a closed oriented
surface of genus G with degree 1 hypersphericity less than 1. 
Then $\Sigma$ contains homologically independent curves $C_1$,
..., $C_G$ with the length of $C_k$ bounded by $200k$.
\end{theorem}

I don't know whether this theorem can be improved to bound the
length of $C_k$ by a constant independent of k and G.  By the
same method, we will prove a bound on the Uryson Width of a
surface of small hypersphericity in a given genus.

\begin{theorem} If $(\Sigma, g)$ is a closed oriented surface of
genus G and degree 1 hypersphericity less than 1, then $(\Sigma,
g)$ has Uryson 1-Width less than 200G + 12.
\end{theorem}

In this theorem, the important parameter is the dependence of the
Uryson 1-Width on the genus.  Given our systolic inequality, it
is rather easy to bound the Uryson 1-Width by a multiple of
$4^G$.  To get the linear bound above, we will have to work
harder.  On the other hand, in section 5 we will give examples of
surfaces with degree 1 hypersphericity less than 1 and Uryson
1-Width on the order of $\sqrt{G}$.  

Both theorems follow from a lemma about surfaces in which each
large ball is a planar domain.  We define the planarity radius of
a surface $(\Sigma, g)$ as the supremal number R so that each
metric ball of radius R is a planar domain.  When the planarity
radius of a surface is large compared to its hypersphericity, an
analogue of Theorem 1.1 holds in any genus.

\begin{lemma} Let $(\Sigma, g)$ be a closed surface with planarity
radius 100 and hypersphericity less than 1.  Then the
Uryson width of $(\Sigma, g)$ is bounded by 50.  Also, there is a
set of disjoint curves on $\Sigma$, each of length less than 24,
whose complement is a planar domain.
\end{lemma}

\proof Since this proof is rather long, I have divided it into
six steps.

Claim 1. Every ball in $\Sigma$ of radius 80 has Uryson 1-Width
less than 12.

pf. Let $\Sigma_1$ be the metric space formed from a ball of
radius 90 around p by collapsing each component of the boundary
to a point.  The manifold $\Sigma_1$ is a topological 2-sphere. 
Any contracting map from $\Sigma_1$ to the unit sphere extends to
a contracting map of the same degree from all of $\Sigma$.  To
see this, let A be the annulus consisting of those points whose
distance from p is between 90 and 100.  The boundary of A
consists of two parts: the boundary of the ball of radius 90,
which we will call B, and the boundary of the ball of radius 100,
which we will call C.  Each component $A_i$ of A is a planar
domain, with one boundary component in B and some number of
boundary components in C.  A contracting map f from $\Sigma_1$ to
the unit sphere is a map from the ball of radius 90, taking each
component $B_i$ of B to a single point $x_i$ in the unit sphere. 
Let $x$ be a point in the unit sphere and $g_i$ a minimal
geodesic from $x_i$ to $x$.  We extend f to A by mapping each
component $A_i$ of the annulus to the geodesic $g_i$,
parametrizing $g_i$ linearly according to the distance from p. 
This map is a contracting map from the ball of radius 100 to the
unit sphere, taking the entire boundary of the ball to x.  We now
extend our map to all of $\Sigma$ by mapping the complement of
$B(p, 100)$ to x.  The new map is contracting and has the same
degree as f.  Therefore, the hypersphericity of $\Sigma_1$ is
less than 1.  By Theorem 1.1, each connected component of a
distance sphere around p in $\Sigma_1$ has diameter less than 12
in $\Sigma_1$.  It follows that each connected component of a
distance sphere around p of radius less than 84 has diameter less
than 12 in $\Sigma$.  This finishes the proof of claim 1.

We define $\Gamma$ to be the subgroup of the fundamental group of
$\Sigma$ generated by lasso-shaped curves consisting of an
arbitrary neck followed by a loop of length less than 24.  It is
easy to check that this subgroup is normal, and we call the
quotient group Q.  Let $\hat \Sigma$ be the cover of $\Sigma$
corresponding to $\Gamma$, and note that Q acts on $\hat \Sigma$
with quotient $\Sigma$.  The fundamental group of $\hat \Sigma$ is
$\Gamma$ and is generated by lassos in $\hat \Sigma$ with loops of
length less than 24.

Claim 2. $\hat \Sigma$ is a planar domain.

pf. Let B be any ball of radius 80 in $\Sigma$.  The fundamental
group of B is generated by lassos with loops equal to the
boundary circles of B.  By Claim 1, each of these boundary
circles has diameter less than 12 in $\Sigma$, so each lasso can
be factored into a composition of lassos with loops of length
less than 24.  Therefore, the inverse image of B in $\hat \Sigma$
consists of isometric copies of B.  Every ball of radius 80 in
$\hat \Sigma$ is an inverse image of a ball of radius 80 in
$\Sigma$, and is a planar domain.  Any two curves in $\hat
\Sigma$ of length less than 24 have intersection number zero,
because if they intersect at all, then they both lie in a ball of
radius 80, and that ball is a planar domain.  On the other hand,
any curve in $\hat \Sigma$ is homologous to a union of curves of
length less than 24.  Therefore, the intersection number of any
two curves in $\hat \Sigma$ is zero.  It follows that $\hat
\Sigma$ is a planar domain.

Claim 3. $\hat \Sigma$ has Uryson width less than 24.

pf. In order to prove this claim, we extend Theorem 1.1 to
complete metrics on planar domains, with a somewhat worse
constant.  In the non-compact case, we define the hypersphericity
in terms of compactly supported maps - maps taking the complement
of a compact set in the domain to a base point on the range. 
(Compactly supported maps have a well-defined degree.)  Our proof
of Theorem 1.1 does not quite extend to this case, because the
map we construct takes the complement of a compact set to the
equator of the target sphere.  We can repair the problem easily
by composing our map with a degree 1 map of the 2-sphere that
takes the equator to a point.  Using the exponential map, we can
construct such a map with Lipshitz constant 2, so we have proved
that a planar domain with hypersphericity 1 has Uryson 1-Width
less than 24.  Moreover, our construction yields a map which is
supported inside a ball of radius 24.  Therefore, we have proven
that if every ball of radius 24 in a planar domain has
hypersphericity less than 1, then the whole surface has Uryson
width less than 24.  Since $\hat \Sigma$ is a planar domain, and
every ball of radius 80 in $\hat \Sigma$ has hypersphericity less
than 1, Claim 3 follows.

Recall from section 2 that a curve in a Riemannian manifold is
called strictly straight if any minimal segment joining two
points on the curve is contained in the curve.  Now let S be the
set of all strictly straight circles in $\hat \Sigma$, of length
no more than 24.  Since $\hat \Sigma$ is planar, the strictly
straight circles in $\hat \Sigma$ are disjoint by Lemma 2.2. (The
union of all the strictly straight circles in S is a closed set,
which we also call S.)

Claim 4. Any curve $\gamma$ from one end of $\hat \Sigma$ to
another intersects a circle in S.

pf. We do this proof in two cases.  First we deal with the
special case that the genus of $\Sigma$ is 1.  By Theorem 4.1,
the shortest homologically non-trivial geodesic in $\Sigma$ has
length less than $4 \pi$.  We pick a point p on this geodesic and
consider the ball of radius 80 around a lift of p.  We see that
$\hat \Sigma$ is diffeomorphic to $R \times S^1$, and that the
strictly straight geodesics include a circle in the class
$[S^1]$.

Second, we deal with the main case that the genus of $\Sigma$ is
at least 2.  We let c be a curve of length less than 24 with
non-zero intersecton number with $\gamma$.  On a compact
manifold, we know that there is a curve homologous to c which
consists of a union of shorter strictly straight circles.  This
is not always true on a non-compact manifold, but we will prove
that it is true on $\hat \Sigma$.  We first homotope the circle c
to a minimal geodesic in its free homotopy class.  Because the
genus of $\Sigma$ is at least 2, the circle c cannot be moved
more than a bounded distance without stretching, and so there is
an actual minimizer.  If this minimizer is not strictly straight,
we surger it along a minimal segment into two circles.  Then we
repeat the procedure, homotoping each circle to a minimal
geodesic in its free homotopy class, and so on.  Because of
bounded geometry, each time we do a surgery, splitting a closed
geodesic of length less than 24 into two pieces, neither piece
null-homotopic, and then homotoping the pieces to minimal
geodesics, the length of each new piece is less than
$(1-\epsilon)$ times the length of the old circle.  Because the
injectivity radius of $\hat \Sigma$ is bounded below, this
procedure must terminate after a finite number of surgeries, with
a union of strictly straight circles each of length less than 24. 
Since the union is homologous to c, one of the curves intersects
$\gamma$.

Claim 5. Q has no torsion.

pf. Suppose Q had a torsion element q of order n.  Pick a point $x$
in $\hat \Sigma$ and let $\gamma$ be a minimal geodesic from $x$ to
$q \cdot x$.  The union of $\gamma$, $q \cdot \gamma$, $q^2 \cdot
\gamma$, ..., and $q^{n-1} \cdot \gamma$ is a circle c in
$\hat \Sigma$.  If we parametrize $\gamma$ by $[0,1]$, then there
is a corresponding parametrization of c by $[0,n]$, so that two
points in c differing by an integer are related by a power of q. 
By Claim 3, there is a map $\pi$ from $\hat \Sigma$ to a tree T,
whose fibers have diameter less than 24.  Applying $\pi$ to the
parametrization of c, we get a map $\pi_0$ from $[0,n]$ to the
tree T.  We will modify this map $\pi_0$ so that each interval
$[k,k+1]$ (for k an integer) meets each vertex of T in a
connected segment, whose image in $\hat \Sigma$ has length no
more than 24.  If the preimage of a vertex t of T intersected
with the interval $[k,k+1]$ consists of several numbers, we
define $a_{min}$ to be the smallest of the numbers and $a_{max}$
to be the largest of the numbers.  Then we modify $\pi_0$ so that
it maps the interval $[a_{min}, a_{max}]$ to the vertex t.  Since
the image of $[k, k+1]$ is a minimal geodesic in $\hat \Sigma$
and the diameter of each fiber of $\pi$ is less than 24, the
length of the image of the interval $[a_{min}, a_{max}]$ is less
than 24.  We repeat this procedure for each vertex t in T.  We
refer to the resulting map from $[0,n]$ to T as $\tilde \pi$.  

The map $\tilde \pi$ obeys two estimates.  First, the inverse
image of any vertex intersected with any interval $[k, k+1]$
consists of a closed sub-interval.  Second, $\tilde \pi$ agrees
with $\pi$, except on the union of finitely many disjoint open
intervals, each of which is mapped to a single vertex of t and
each of which parametrizes a segment of length less than 24. 
Now, it is a topological fact that two numbers in $[0,n)$
differing by an integer are mapped to the same point of T by
$\tilde \pi$.  Call the two numbers a and b, and let p and q be
the points parametrized by a and b.  It may not be the case that
p and q are mapped to the same point by $\pi$, because we have
modified $\pi$ along some disjoint intervals.  Let $a_0$ and
$b_0$ be the closest endpoints of the intervals containing a and
b, and let $p_0$ and $q_0$ be the points in c parametrized by
$a_0$ and $b_0$.  The distance from p to $p_0$ is less than 12,
as is the distance from q to $q_0$.  Also, $p_0$ and $q_0$ are
mapped to the same point of T by the original $\pi$, so the
distance between them is no more than 24.  Therefore, the
distance between p and q is no more than 48.  On the other hand,
the distance between any two different points in the same Q-orbit
is at least 80.  This contradiction establishes the claim. 

Claim 6. Any two points in the same Q orbit are in different
components of $\hat \Sigma - S$.

pf. Let $\gamma$ be any curve connecting $x$ to $q \cdot
x$, avoiding all the circles in S.  Since there is no torsion in
Q, the union of all the curves $q^n \cdot \gamma$, for all
integers n, is a single curve running from one end of $\hat
\Sigma$ to another.  By Claim 4, this curve must intersect a
circle in S.  But since S is invariant under the action of Q (or
under any isometry), $\gamma$ must intersect a circle in S.

Because S is Q-invariant, the curves in S descend to disjoint
curves in $\Sigma$ of length less than 24, dividing $\Sigma$ into
planar domains.  This establishes the last statement of the
lemma.  We now bound the Uryson width of $\Sigma$.  We select
some union of components of $\hat \Sigma - S$ that form a
fundamental domain for the action of Q on $\hat
\Sigma - S$, and on each one we construct a map to a tree for
which each boundary circle is a fiber, and with fibers of
diameter less than 50.  On each component, we pick a point p and
begin by taking the distance spheres around p as our fibers (or
more properly, their intersection with the component we have
selected).  These have diameter less than 24, but they do not
extend to the boundary with each boundary component as a fiber. 
We surger the fibers near each boundary circle, so that the
boundary circle is a fiber.  This increases their diameter by no
more than 24, and so the Uryson Width of $\Sigma$ is less than
50. \endproof

Now we turn to the applications of Lemma 4.4.  If a surface has
planarity radius less than 100, then it contains two smooth
embedded curves, each of length less than 200, intersecting
transversely at one point.  To see this, we find a non-planar
ball of radius less than 100 in our surface.  We slightly change
the metric on the ball to make the boundary convex and then
consider a minimal rational basis of the first homology group of
the ball.  Each curve in the ball is homologous to a sum of
curves of length less than 200, and so every curve in the minimal
basis will have length less than 200.  According to Lemma 2.2,
each curve will be strictly straight, and so any pair of curves
will intersect no more than once, and if they do intersect once
the intersection will be transverse.  Since the domain is not
planar, two curves in the basis must intersect, and this proves
the claim.

Therefore, an arbitrary surface has either two short embedded
curves intersecting transversally at one point or a large
planarity radius.  Exploiting these two structures, we can quickly
prove Theorems 4.2 and 4.3.

\vskip5pt

Proof of Theorem 4.2:  If the planarity radius of $\Sigma$ is
less than 100, then there are two curves of length less than 200
with intersection number one, intersecting in a single point.  We
call these curves $C_1$ and $C_2$.  Removing these curves from
$\Sigma$ leaves a manifold U whose boundary is a single circle
(with corners).  We can embed U into a closed manifold $\Sigma'$
of genus G-1, so that the boundary circle bounds a disk in
$\Sigma'$.  Now we extend the metric on U to give a Riemannian
metric on $\Sigma'$, but make this metric very small on the
complement of U.  If we make the metric sufficiently small on the
complement of U, then there is a degree 1 map from $\Sigma$ to
$\Sigma'$ with Lipshitz constant arbitrarily close to 1, and so
the degree 1 hypersphericity of $\Sigma'$ is less than 1.

By induction on the genus, $\Sigma'$ contains $G-1$ homologically
independent curves $C_1'$, ... $C_{G-1}'$ with the length of
$C_k'$ bounded by $200k$.  By taking a minimal basis of homology
of $\Sigma'$, we can assume that each curve is strictly straight,
and by altering them very slightly we can assume that each curve
passes through the complement of U at most once.  Now we define a
curve $C_{k+2}$ in $\Sigma$ as follows.  Begin with the portion
of the curve $C_k'$ lying in U and lift it to $\Sigma$.  If this
curve is closed, we are done.  If it is not closed, then it is a
single segment with endpoints on the union of $C_1$ and $C_2$. 
Since the diameter of this union is 200, we can close the curve
by adding a portion of $C_1$ and a portion of $C_2$ of total
length no more than 200.  Therefore, the length of $C_{k+2}$ is
less than $200 (k+1)$.  We have now defined curves $C_1$, ...
$C_{G+1}$ in $\Sigma$ with the length of $C_k$ bounded by $200k$.

We now show that the curves $C_k$ are homologically independent. 
Let S be the rational span of their homology classes in
$H_1(\Sigma)$.  The map from $\Sigma$ to $\Sigma'$ induces a map
from S to $H_1(\Sigma')$.  The kernel of this map is at least
two-dimensional, generated by $C_1$ and $C_2$.  On the other
hand, the curve $C_{k+2}$ is mapped to a curve homologous to
$C_k'$, and so the image of this map has dimension at least G-1. 
Therefore, the space S has dimension G+1 and the curves $C_k$ are
homologically independent.

On the other hand, if the planarity radius of $(\Sigma, g)$ is
more than 100, then the lemma guarantees G homologically
independent curves of length less than 24.  \endproof

Proof of Theorem 4.3: If the planarity radius of $(\Sigma, g)$ is
less than 100, then we again find two curves of length less than
200 intersecting transversally in a point, define U to be the
manifold with boundary formed by removing these two curves, embed
U in a surface $\Sigma'$ of genus G-1, and extend the metric on U
to a metric on $\Sigma'$ which is very small away from U.  We can
then define a degree 1 map F from $\Sigma$ to $\Sigma'$ with
Lipshitz constant only slightly more than 1, and with the
property that the preimage of a set of diameter D has diameter
less than $D+200$.  (The map F is essentially the map which
collapses the two curves to a point.)  The degree 1
hypersphericity of $\Sigma'$ is less than 1, and so by induction
on the genus we may assume that there is a map $\pi'$ from
$\Sigma'$ to a graph whose fibers have diameter less than $200
(G-1) + 12$.  (The base for the induction, G=0, is just Theorem
1.1.)  We define $\pi$ to be the composition of $\pi'$ with F.  The
diameters of the the fibers of $\pi$ are less than 200G + 12, and
in this case the theorem holds.  On the other hand, if $(\Sigma,
g)$ has planarity radius more than 100, then by the lemma it has
Uryson 1-width less than 50, and the theorem holds in this case
too. \endproof

\section{Strange Metrics on Surfaces of High Genus}

In this section we will construct surfaces of arbitrarily small
hypersphericity and with Uryson 1-Width at least 1.

\begin{theorem} There exists a sequence of closed oriented
surfaces $(\Sigma_i, g_i)$ with hypersphericity tending to 0 and
Uryson 1-Width at least 1.
\end{theorem}

The construction is roughly as follows.  First, we pick a space
$(M,g)$ along with a non-zero homology class h in $H_2(M,\mathbb{Z})$,
which can be realized by an oriented surface.  We pick $\Sigma_1$
to be any surface in M in the homology class h.  The surface
$\Sigma_i$ is $\Sigma_1$ joined to a very dense mesh of thin
tubes filling M, and the metric $g_i$ is the induced metric from
M.  The resulting surfaces can be chosen to converge to M in the
Gromov-Hausdorff metric.  The details of this construction are
carried out by Cassorla in \cite{Cass}.

Using the geometry of M, we will prove inequalities about the
hypersphericity and Uryson 1-Width of the surfaces $\Sigma_i$. 
Our main tool to relate the geometry of $\Sigma_i$ to the
geometry of M will be Gromov's extension theorem for maps of
small width, and some minor generalizations.  Recall that the
width of a map is defined to be the supremal diameter of its
fibers.  We now recall Gromov's extension theorem for maps of
small width, which appeared in \cite{G1}.

\begin{reftheorem} [Gromov] If p is a surjective map from X to Y
with width less than $\pi/2$, and if f is a contracting map from
X to a unit sphere (of any dimension), then f is homotopic to a
map that factors through p.
\end{reftheorem}

\begin{refcor}[Gromov] If $(M^n, g)$ has Uryson (k-1) width less
than $\pi/2$, then every contracting map from M to the unit k-sphere
is null-homotopic.  
\end{refcor}

As a special case of the corollary, we see that the
hypersphericity of a Riemannian n-manifold is bounded by $(\pi/2)$
times its Uryson (n-1)-Width, which is an inequality we used
several times in this paper.

Using Gromov's extension theorem, we can construct surfaces with
arbitrarily small degree 1 hypersphericity but with
hypersphericity at least 1.  For example, let $(\Sigma, g)$ be a
degree 2 branched covering of the unit 2-sphere, branched over an
$\epsilon$-dense set of points, with the induced metric (strictly
speaking the induced metric is singular, but we can add a tiny
positive symmetric form to make a Riemannian metric).  The
branched covering map has width less than $3 \epsilon$, and it's
certainly surjective.  By Gromov's theorem, any contracting map
from $\Sigma$ to the sphere of radius $3 \epsilon$ factors
through the branched covering and has even degree.  Therefore,
the degree 1 hypersphericity of $\Sigma$ is less than $3
\epsilon$.  On the other hand, the branched covering is itself a
contracting degree 2 map to the unit sphere, and so the
hypersphericity of $\Sigma$ is at least 1.

We will define a class of maps called maps of small girth which
generalize the surjective maps of small width, and extend
Gromov's theorem to the maps of small girth.  This generalization
is not very serious, but it will be convenient for our purposes,
and we will include more details of the proof than appear in
\cite{G1}.  To arrive at the definition of the girth of a map, we
make a slight shift in point of view.  Instead of a polyhedron Y, the
target of our map is a polyhedron Y equipped with a cover by open
sets, $Y = \cup U_i$.  Instead of taking a surjective map from X
to Y, we take a map p whose image intersects each open set $U_i$ in
the cover.  For any point y in Y, let $U_y$ be the union of all
$U_i$ containing y.  Instead of considering the fiber
$p^{-1}(y)$, we consider the fattened fiber $p^{-1}(U_y)$.  We
say that the map p has girth less than w if $p^{-1}(U_y)$ has
radius less than w.

We will say that a Riemannian manifold M has convexity radius at
least w if each ball of radius at most w is convex with respect
to minimal geodesics.  A set U in M is convex with respect to
minimal geodesics if for any two points x and y in U, there is a
unique minimal geodesic in M joining x to y, and this minimal
geodesic lies in U.  For instance, any complete Riemannian
manifold with sectional curvature less than 1 and injectivity
radius at least $\pi$ has convexity radius at least $\pi/2$. 

Now we can prove one version of the extension lemma.

\begin{lemma} If p is a map from X to Y with girth less than w, f
is a contracting map from X to M, and M is a Riemannian manifold
with convexity radius at least w, then there is a map g from Y to
M so that f is homotopic to $g \circ p$.  Moreover, the homotopy
moves each point of X by no more than 2w.
\end{lemma}

\proof First we replace the $U_i$ by a locally finite sub-cover,
and we redefine $U_y$ accordingly.  Since the new $U_y$ is
contained in the old $U_y$, the fattened fibers $p^{-1}(U_y)$
still have radius strictly less than w.  We pick a number
$\epsilon$ so that the girth of p is less than $w - \epsilon$.

We define a map G from Y to the set of convex sets in M.  For
each point y in Y, $G(y)$ is defined to be the convex hull of the
$\epsilon$-neighborhood of $f(p^{-1}(U(y)))$.  By the definition
of girth, $p^{-1}(U(y))$ lies inside a ball of radius $w -
\epsilon$ in X.  Since f is contracting, $f(p^{-1}(U(y))$ lies
inside a ball of radius $w - \epsilon$ in M, and its $\epsilon$
neighborhood lies in a ball of radius w.  Since each ball of
radius w in M is convex, G(y) is contained in a ball of radius w. 
The purpose of the $\epsilon$ neighborhood is to ensure that G(y)
is open.

The main point of the proof is to construct a ``section'' of the
map G: that is, a continuous map g from Y to M so that $g(y)$
lies in $G(y)$ for each y.  This step is a bit technical.  I
think the moral reason why we can find a section is just that
each set G(y) is contractible, but there are technicalities
because the sets $G(y)$ may not vary continuously with y. 

We divide Y into simplices so small that each simplex lies in one
of the $U_i$.  Then we construct the map g one skeleton at a
time, so that at each stage it is continuous and g(y) lies in
G(y).  The map from the zero skeleton is trivial to construct. 
By induction, it suffices to extend our section from the boundary
of a simplex to its interior.  We pick a simplex $\Delta$ and
assume that g is defined on its boundary so that $g(y)$ lies in
$G(y)$ for each y on the boundary.  We choose a homeomorphism of
$\Delta$ with the ball of radius 1.  Since the simplex lies in
one of the $U_i$, there is a point m in M which lies in $G(y)$
for every y in $\Delta$.  We define g(y) to be m for each y in a
ball of radius $1 - \delta$, for some small number $\delta$ which
we will choose later.  Then, we extend g to the annulus $S^{k-1}
\times (1-\delta, 1)$ by mapping each ray $\theta \times (1-\delta,
1)$ to the minimal geodesic between $g(\theta,1)$ and m.  Our extension
g is clearly continuous, and it only remains to check that g(y)
lies in G(y) for every y.

The boundary of $\Delta$ is compact, so it meets only finitely
many of the $U_i$.  We number these $U_1, ..., U_n$.  For any
subset i of the numbers from 1 to n, we define $U_I$ to be the
set of points in the boundary of the simplex which lie in $U_i$
exactly when i is a member of I.  We define $G_I$ to be the
convex hull of the epsilon-neighborhood of $f(p^{-1}(\cup_{i \in
I} U_i))$.  The condition that $g(y)$ lies in $G(y)$ means
exactly that g maps $U_I$ into $G_I$ for each I.  Now the sets
$U_I$ are not closed, but each point in the closure of $U_I$ lies
in a $U_J$ where J is a subset of I, and so we see that g maps
the closure of $U_I$ into $G_I$ for each I.  Since the sets $G_I$
are open, g maps a small neighborhood of the closure of each
$U_I$ into $G_I$.  Therefore, we may choose closed sets $V_i
\subset U_i$, so that when we define $V_I$ as the set of points y
which lie in $V_i$ exactly when i is a member of I, g(y) lies in
$G_I$ for each y in $V_I$.  This tiny improvement gives us the
space we need to prove that g(y) lies in G(y).  We define $W_i
\subset S^{n-1}$ to be the intersection of $V_i$ with $S^{n-1}$,
and we choose a number $\delta$ so small that $W_i \times [0,
\delta]$ lies in $U_i$ for each i.  For each point $\theta$ in
$S^{n-1}$ and each t in $(1-\delta, 1)$, $V_{(\theta,1)}$ is a
subset of $U_{(\theta, t)}$ and $g(\theta,1)$ lies in
$G(\theta,t)$.  Since m also lies in $G(\theta,t)$ for each t, we
see that the minimal geodesic from m to $g(\theta,1)$ lies in
$G(\theta,t)$, and therefore g(y) lies in G(y) for each y.

Both $f(x)$ and $g(p(x))$ lie in $G(p(x))$ which lies in a ball
of radius w.  Since the convexity radius of M is at least w,
$f(x)$ and $g(p(x))$ are joined by a unique minimal geodesic, and
this minimal geodesic varies continuously with x.  The family of
geodesics, parametrized proportionally to length, is a homotopy
from f to $g \circ p$.  This homotopy moves each point of X less
than 2w.
\endproof

With this version of the extension lemma, we can prove estimates
relating the Uryson 1-widths and hypersphericities of the
surfaces $\Sigma_i$ to the geometry of M.

We prove that the Uryson 1-Width of each surface $\Sigma_i$ is at
least the convexity radius of M.  To see this, suppose that
$\Sigma_i$ admitted a map p to a 1-polyhedron Y whose fibers had
diameter less than the convexity radius of M.  After changing the
target, we can assume that p is surjective.  Therefore, the girth
of p is less than the convexity radius of M.  The inclusion of
$\Sigma_i$ in M is a contracting map (because $\Sigma_i$ has the
induced metric from M).  Therefore, the extension lemma applies,
proving that the inclusion of $\Sigma_i$ in M factors through p
up to homotopy.  But the map p kills the fundamental homology
class of $\Sigma_i$, and the inclusion of $\Sigma_i$ does not. 
Therefore, the Uryson 1-Width of $\Sigma_i$ is at least the
convexity radius of M.

Next we bound the hypersphericity of $\Sigma_i$ for large values
of i.  Because the surfaces $\Sigma_i$ are converging to M in the
Gromov-Hausdorff metric, for large values of i, the inclusion of
$\Sigma_i$ in M has girth $\epsilon_i$ tending to 0.  Now, by the
extension lemma, any contracting map from $\Sigma_i$ to the
sphere of radius $\epsilon_i$ factors up to homotopy through the
inclusion of $\Sigma_i$ in M.  If every map from M to the
2-sphere kills the homology class h, then the hypersphericity of
$\Sigma_i$ is less than $\epsilon_i$.  For example, if M is the
complex projective plane and h is any non-zero homology class in
$H_2(M,\mathbb{Z})$, then every map from M to the 2-sphere kills h.

If we take M to be the complex projective plane with its standard
metric, these two estimates prove Theorem 5.1.

The constructions that we have made for surfaces of high genus
also apply to manifolds of higher dimension.  The
Gromov-Hausdorff convergance is provided by a theorem of Ferry
and Okun from \cite{FO}.  We state a weak version of their
theorem, which suffices for our purposes.

\begin{reftheorem} [Ferry and Okun] Let A be a connected manifold of
dimension at least 3 and B be a connected Riemannian manifold,
and let f be a map from A to B which is surjective on $\pi_1$. 
Then there is a sequence of metrics $g_i$ on A and of maps $f_i$
from A to B homotopic to f, so that $(A, g_i)$ converges in the
Gromov-Hausdorff metric to B by the maps $f_i$.
\end{reftheorem}

The construction is roughly as follows.  First, homotope f so
that the image contains a very dense set of points in B.  Then,
homotope f along various curves joining the preimages of this
dense set, so that each nearby pair of points in the dense set is
joined by a minimal geodesic lying in the image of f.  This last
step requires us to find curves in the domain which map to the
right homotopy class in the range, which requires surjectivity on
$\pi_1$.

In particular, this theorem implies that the maps $f_i$ have
girth tending to 0 and Lipshitz constant tending to 1.  Using
this theorem, we can construct metrics on the 4-sphere with
arbitrarily small hypersphericity and Uryson 3-width at least 1. 
To see this, take A to be the 4-sphere, B to the quaternionic
projective plane, and f to be the hyperplane embedding.  Applying
the theorem of Ferry and Okun, we get metrics $g_i$ on A and maps
$f_i$ with girth tending to 0.  Applying the estimates just as
for surfaces, we conclude that the Uryson 3-width of $(S^4, g_i)$
is at least the convexity radius of the quaternionic projective
plane and that its hypersphericity tends to 0.

We will give two other sequences of surfaces with hypersphericity
tending to 0 and Uryson 1-width at least 1.  These examples are
included to show that estimating the hypersphericity of a surface
is rather complicated.  In each example, a different piece of
geometry or topology is used to rule out a non-zero degree
contracting map to a small sphere.  In our first example, the
obstruction came from homotopy theory.  In the second example, it
will come from lattice theory.  In the third example, it will
come from some considerations related to the theory of systoles. 
These examples require more work to construct then our first
example, but they also have slightly stronger properties.  For
instance, we will construct surfaces $\Sigma_i$ so that the
smallest Lipshitz constant of a homotopically non-trivial map to
$CP^{\infty}$ tends to infinity while the Uryson 1-Width remains
bounded below.

The estimates in these examples depend on a modification of the
extension lemma, giving a Lipshitz extension of a Lipshitz
function.  In order to state the modified lemma, we make some
definitions.  We say that a map p is dense on the scale S if the
image of p meets every ball of radius S.  We say that a map p is
expanding on the scale S if, for any two points $x_1$ and $x_2$
in the domain, $dist(p(x_1), p(x_2)) > dist(x_1, x_2) - S$.  For
points in the domain closer together than S, this inequality is
vacuous, but for points in the domain much farther apart than S,
the map p is nearly expanding.  We say that a Riemannian manifold
has bounded geometry on a scale S if its sectional curvature is
pinched between $1/S^2$ and $-1/S^2$, and if its injectivity
radius is at least $\pi S$.  For instance, the unit sphere has
bounded geometry on the scale 1.

\begin{lemma} Let X and Y be locally compact path metric spaces. 
Let p be a map from X to Y which is dense on the scale W and
expanding on the scale W, for some number W less than $1/200$. 
Let M be a complete Riemannian n-manifold with bounded geometry
on the scale 1.  Let f be a contracting map from X to M.  Then
there is a map h from Y to M, with Lipshitz constant $50n$, so
that f is homotopic to $h \circ p$.
\end{lemma}

\proof First we prove that the girth of p is bounded by 4W.  To
see this, we let $U_i$ be the cover of Y by W-balls.  The set
$U_y$ is defined as the union of all W-balls meeting y, which is
the 2W-ball centered at y.  Since the image of p is W-dense, the
set $U_y$ is contained in a 3W-ball around a point in the image
of p, say p(x).  Since p is expanding on the scale W, the
preimage of $U_y$ is contained in the 4W ball around x.  Because
the girth of p is bounded by 4W, which is less than 1/50, and the
convexity radius of M is at least $\pi/2$, we can apply Lemma
5.2, which constructs a map g from Y to M obeying the estimate
$dist(f(x), g(p(x)) < 8W$.

We will prove that the map g which we constructed in Lemma 5.2
obeys a quasi-Lipshitz inequality, bounding the distance between
$g(y_1)$ and $g(y_2)$ by $dist(y_1, y_2) + 19W$.  Since p is
dense on the scale W, there are points $p(x_1)$ and $p(x_2)$
within W of $y_1$ and $y_2$, respectively.  Recall from the proof
of Lemma 5.2 that the set G(y) is the convex hull of a small
neighborhood of $f(p^{-1}(U_y))$, and that g(y) lies in G(y). 
Since $p(x_1)$ lies in $U_{y_1}$, the point $f(x_1)$ lies in
$G(y_1)$.  Since the radius of G(y) is less than 4W, the distance
between $f(x_1)$ and $g(y_1)$ is less than 8W.  The distance from
$g(y_1)$ to $g(y_2)$ is less than $dist(g(y_1),f(x_1)) +
dist(f(x_1),f(x_2)) + dist(f(x_2),g(y_2))$, which is less than
$dist(f(x_1),f(x_2)) + 16 W$.  Since f is contracting, the
distance from $f(x_1)$ to $f(x_2)$ is less than the distance from
$x_1$ to $x_2$.  Since p is expanding on the scale W, this
distance is less than $dist(p(x_1),p(x_2)) + W$.  But the point
$p(x_1)$ lies within W of $y_1$, and the point $p(x_2)$ lies
within W of $y_2$.  Putting all these estimates together, we get
the quasi-Lipshitz inequality.

We now explain how to build a Lipshitz map h on the scaffold of
the quasi-Lipshitz map g.  This construction is adapted from the
paper \cite{LPS} of Lang, Pavlovic, and Schroeder, in which they
construct Lipshitz extensions of Lipshitz maps into negatively
curved spaces.  For each y in Y, we define a set $H(y)$ in M as
the intersection of the closed balls $B(g(y'), 2 dist(y,y') +
20W)$ for every y' in Y.  Because of our quasi-Lipshitz
inequality, the ball of radius W around g(y) lies in H(y), and so
H(y) is not empty.  On the other hand, H(y) lies within the ball
of radius 20W around g(y), a ball of radius less than 1/10.  We
can think of H(y) as showing the rough location of h(y).  Notice
that if the distance from y to y' is at least 20W, then the ball
$B(g(y'), 2 dist(y,y') + 20 W)$ contains the ball around g(y) of
radius 20W, which in turn contains H(y).  Therefore, H(y) is an
intersection of balls of radius at most 60W.  Because of the
bounded geometry of M, each of these balls is convex, and so
$H(y)$ is a convex set.

Unlike the sets G(y), the sets H(y) vary nicely with y.  It is
easy to see that H(y) varies continuously with y.  Because each
set H(y) contains a ball of radius W, H(y) is a convex set with
interior.  We define the function d(m) to be the infimum of the
expression $2 dist(y, y') + 20 W - dist(m, g(y'))$ as y' varies
over Y.  The function d(m) is continuous.  It is positive on the
interior of H(y), zero on the boundary, and negative on the
complement of H(y).  The set H($y_0$) contains the set where $d(m)$
is greater than $2 dist(y,y_0)$ and is contained in the set where
$d(m)$ is greater than $- 2 dist(y,y_0)$.  Since d is continuous,
we see that H(y) varies continuously with y, using the Hausdorff
topology on closed subsets of M.  We will show something much
stronger, namely that the Hausdorff distance between $H(y)$ and
$H(y')$ is bounded by $5 dist(y,y')$.  This proof is very closely
modeled on a similar proof in \cite{LPS}, in a slightly different
situation.  Let m' be a point in H(y') and let m be the closest
point to m' in H(y).  We will bound the distance from m to m' by
5 dist(y,y').  Since the argument applies to each point m' in
H(y') and by symmetry to each point in H(y), this estimate will
bound the Hausdorff distance between H(y) and H(y') by 5
dist(y,y').

Since m lies on the edge of $H(y)$, it must lie on the edge of
some of the balls defining $H(y)$.  We will call a point $y_1$ in
Y taut if m lies on the edge of $B(g(y_1), 2 dist(y,y_1) + 20W)$. 
If $y_1$ is a taut point, then we will call the ray from m to
$g(y_1)$ a taut ray.  We proved above that if the distance from
$y_1$ to y is more than 20W, then H(y) lies inside $B(g(y_1), 2
dist(y, y_1) + 20W)$.  Therefore, every taut point lies within
20W of y, and every taut ray has length less than 60W, which is
less than 1/3.  We will prove that the angle between any two taut
rays is no more than $(.34) \pi$.  

Let $y_1$ and $y_2$ be two taut points.  We define the vectors a
and b to be the inverse images of $g(y_1)$ and $g(y_2)$ under the
exponential map of M at m.  The absolute value of a equals the
distance from m to $g(y_1)$, which equals $2 dist(y_1, y) + 20W$
by the definition of a taut ray.  Similarly, the absolute value
of b equals the distance from m to $g(y_2)$, which equals $2
dist(y_2, y) + 20 W$.  By the quasi-lipshitz estimate for g, the
distance from $g(y_1)$ to $g(y_2)$ is less than $dist(y_1,y_2) +
20$, which is less than half of $|a| + |b|$.  Since taut rays
have length less than 1/3, $g(y_1)$ and $g(y_2)$ lie inside the
ball of radius 1/3 around m, and since this ball is convex, the
minimal geodesic between them lies in the ball as well.  By the
Rauch comparison theorem, within the ball of radius 1/3, the
exponential map does not contract lengths by more than the factor
$\frac{1/3}{sin(1/3)}$.  Therefore, the distance from a to b is
no more than $\frac{1/3}{sin(1/3)}$ times the distance from
$g(y_1)$ to $g(y_2)$, and $|a - b| < (54/53) dist(g(y_1), g(y_2))
< (27/53) (|a| + |b|)$.  By trigonometry, the angle between a and
b is less than $(.34) \pi$. 

Because m is the closest point to m' in $H(y)$, the taut rays
from m must contain a ray pointing directly away from m' in their
convex hull.  Since the taut rays cluster together, none of them
can lie too far from this ray pointing directly away from m'.  In
particular, the angle between a taut ray and the ray from m to m'
must be greater than $(.66) \pi$.

It suffices to prove our Lipshitz inequality for y very close to
y', and since H(y) varies continuously with y, we can assume that
m is very close to m'.  Let $y_1$ be a taut point in Y.  Because
m' is very close to m, the distance between them is bounded by
$(1/cos (.34 \pi)) (dist(g(y_1), m') - dist(g(y_1), m))$.  But by
the definition of a taut point, $dist(g(y_1),m) = 2 dist(y_1,y) +
20W$.  On the other hand, since m' lies in H(y'),
$dist(g(y_1),m') \le 2 dist(y_1,y') + 20W$.  Therefore, the
distance from m to m' is less than $(2/cos(.34
\pi)) dist(y,y')$, which is less than 5 dist(y,y').

We now define h(y) to be the center of mass of the
1/10-neighborhood of H(y).  Centers of mass for measures in
Riemannian manifolds are explained in section 4 of \cite{LPS}. 
We prove some estimates very closely following theirs, but in a
slightly different situation.  First, we define center of mass. 
For a probability measure $\mu$ on a Riemannian manifold, we
define a function $D_{\mu}(p) = \int dist(p,q)^2 d\mu(q)$.  If
the manifold has sectional curvature less than 1 and injectivity
radius more than $\pi$, and if the distance from p to q is less
than 2/5, then the function $d(\cdot, q)^2$ is strictly convex at p,
with Hessian greater than 1.  Therefore, if p lies in a ball of
radius 1/5 and $\mu$ is supported in the same ball, $D_{\mu}$ is
also strictly convex at p with Hessian greater than 1.  Since the
ball of radius 1/5 is convex with respect to minimal geodesics,
the gradient of the function $D_{\mu}$ points out at every point
on the edge of the ball.  Therefore, the function $D_{\mu}$ has a
unique minimum in the ball of radius 1/5.  The point where this
minimum is attained is called the center of mass of $\mu$, and we
will denote it $c_{\mu}$.

Because the Hessian of $D_\mu$ is greater than 1 on the ball of
radius 1/5, $D_{\mu}(p) > D_{\mu}(c_{\mu}) + (1/2) dist(p,
c_{\mu})^2$ for each p in the ball.  Using this estimate, we can
bound the distance between the centers of mass of two measures,
supported on the same ball of radius 1/5.

$$dist(c_{\mu}, c_{\mu'})^2 < D_{\mu}(c_{\mu'}) -
D_{\mu}(c_{\mu}) + D_{\mu'}(c_{\mu}) - D_{\mu'}(c_{\mu'}).$$

Plugging in the definition of $D_{\mu}$ we get the following
expression.

$$ = \int (dist(c_{\mu},q)^2 - dist(c_{\mu'},q)^2)(d\mu' -
d\mu).$$

Factoring the difference of squares, we can bound this
expression.

$$< (4/5) dist(c_{\mu}, c_{\mu'}) \int |d\mu' - d\mu|.$$

Canceling one factor of $dist(c_{\mu}, c_{\mu'})$ from each side
leaves us with the following bound.

$$dist(c_{\mu}, c_{\mu'}) < (4/5) \int |d\mu' - d\mu|.$$

The probability measures we will use are the renormalized volume
measures of the 1/10 neighborhood of H(y), which we call
$\mu(y)$.  If y' is sufficiently close to y, then $\mu(y)$ and
$\mu(y')$ will be supported on a ball of radius 1/5.  Again, if
they are sufficiently close, then $\int |d\mu(y') - d\mu(y)|$
will be bounded by the surface area to volume ratio of the 1/10
neighborhood of H(y) times the Hausdorff distance between the
1/10 neighborhood of H(y) and the 1/10 neighborhood of H(y'). 
The surface area to volume ratio of a 1/10 neighborhood of any
set in an n-manifold with $Ric > -(n-1)$ is at most that of the
1/10 ball in hyperbolic n-space, which is less than 11n. 
Finally, the Hausdorff distance between the 1/10 neighborhoods of
H(y) and H(y') is no more than the Hausdorff distance between
H(y) and H(y'), which is less than $5 dist(y,y')$.  Putting the
bounds together, we see that the Lipshitz constant of h is
less than 44n. \endproof

We define the hypersphericity of a homology class h in $H_2(M)$
for a Riemannian manifold $(M,g)$ as the supremal R so that there
is a contracting map from $(M,g)$ to the 2-sphere of radius R which
does not kill the homology class h.  If M is a closed oriented
surface and h is the fundamental class of M, then we recover the
usual definition of hypersphericity.  By Lemma 5.3, if the
hypersphericity of h in M is less than $\epsilon$, then for
sufficiently large i, the hypersphericity of $\Sigma_i$ is less
than $100 \epsilon$.

We will now construct some sequences of higher dimensional
manifolds $(M_n,g_n)$ with bounded geometry on the scale 1 and
with hypersphericity of $h_n$ in $M_n$ tending to 0, where $h_n$
is non-zero and can be realized by an oriented surface.  For each
manifold $M_n$, we have constructed a sequence of surfaces
$\Sigma_i$ in the homology class $h_n$ converging to $(M_n,g_n)$. 
By Gromov's extension lemma, we know that the Uryson 1-Width of
$\Sigma_i$ is greater than 1 for every i.  On the other hand, by
the Lipshitz extension lemma, we know that for sufficiently large
n and sufficiently large i, the hypersphericity of $\Sigma_i$ is
as small as we like.

Our first examples of manifolds with bounded geometry carrying
homology classes of small hypersphericity are high-dimensional
flat tori.  Minkowski discovered that in high dimensions a random
flat torus behaves differently from a rectangular torus.  A
rectangular torus of volume 1 has injectivity radius at most a
half, but a random flat N-torus of volume 1 has injectivity
radius on the order of $\sqrt{N}$.  (See the Minkowski-Hlawka
theorem in a book on the geometry of numbers, such as \cite{Ca}.) 
One reason for the difference in behavior is that in high
dimensions a ball of volume 1 has radius on the order of
$\sqrt{N}$.  A random flat torus usually has a fundamental domain
that looks much more like a ball than like a rectangle.

Let T be a flat torus with volume 1 and injectivity radius on the
order of $\sqrt{N}$.  Any homologically non-trivial curve
in T has length on the order of $\sqrt{N}$.  On the other hand,
the small volume of T puts some restrictions on the number of
contracting maps that T can have to circles of large radii.

\begin{lemma} Let T be a flat N-torus with injectivity radius 1 and
volume less than $(c \sqrt{N})^{-N}$.  Then T has a one dimensional
homology class h, so that any map to the unit circle which does
not kill h has Lipshitz constant at least $2 \pi c \sqrt{N}$.
\end{lemma}

\proof Suppose the conclusion does not hold.  Let $h_1$ be any
1-dimensional homology class, and let $f_1$ be a map to $S^1$ not
killing $h_1$, with Lipshitz constant less than $2 \pi c
\sqrt{N}$.  Then let $h_2$ be a one dimensional homology class
killed by $f_1$, and let $f_2$ be a map to $S^1$ not killing
$h_2$, with Lipshitz constant less than $2 \pi c \sqrt{N}$.  Let
$h_3$ be a 1-dimensional homology class killed by $f_1$ and
$f_2$, and let $f_3$ be a map to $S^1$ not killing $h_3$, with
Lipshitz constant less than $2 \pi c \sqrt{N}$.  Continuing in
this way, we produce a map F of non-zero degree from T to the
direct product of N circles.  Since each coordinate of F has
Lipshitz constant less than $2 \pi c \sqrt{N}$, the volume
dilation of F is less than $(2 \pi)^N (c \sqrt{N})^N$, and the
volume of the image of F is less than $(2 \pi)^N$.  Therefore F
is not surjective and has degree 0. \endproof

Next, we extend our lemma to homology classes of dimension 2 or
higher.

\begin{lemma} Let $T_0$ be a flat N-torus with injectivity radius
1 and volume less than $(c \sqrt{N})^{-N}$, let $T_1$ be the unit
cube (n-1)-torus, and let T be their Cartesian product.  Then T
has an n-dimensional homology class h so that any map from T to
the unit n-sphere which does not kill h has n-dilation at least
$\omega_n c \sqrt{N}$, where $\omega_n$ is the volume of the unit
n-sphere. In particular, the map has Lipshitz constant at least
$(w_n c \sqrt{N})^{1/n}$.
\end{lemma}

\proof Again, we suppose the conclusion is false.  Let a be the
fundamental homology class of $T_1$ in T.  Let $h_1$ be a one
dimensional homology class in $T_0$, and let $f_1$ be a map from
T to the unit n-sphere not killing $h_1 \times a$, with
n-dilation less than $\omega_n c \sqrt{N}$.  Pick $h_2$, a
one-dimensional homology class in $T_0$, so that $f_1$ kills $h_2
\times a$, and let $f_2$ be a map from T to the unit n-sphere not
killing $h_2 \times a$, with n-dilation less than $\omega_n c
\sqrt{N}$.  Proceeding in this way, construct $f_i$ for i from 1
to N.  Now let $T'$ be $T_0 \times T_1^N$, and let $\pi_i$ be the
projection of $T'$ onto a copy of T provided by the product of
$T_0$ with the $i^{th}$ copy of $T_1$.  Now we define a map F
from $T'$ to the Cartesian product of N unit n-spheres.  The map
to the $i^{th}$ unit n-sphere is given by $f_i \circ \pi_i$.  It
is easy to check that F has non-zero degree.  On the other hand,
since each component of F has n-dilation less than $\omega_n c
\sqrt{N}$, the volume dilation of F is less than
$(\omega_n)^N (c \sqrt{N})^N$, and the volume of the image of F
is less than $\omega_n^N$.  Therefore F is not surjective and has
degree zero.  This contradiction finishes the proof of the lemma. 
\endproof

Taking n to be 2, Lemma 5.5 gives us a sequence of flat tori
$T^N$ with injectivity radius 1 and 2-dimensional homology
classes $h_N$, with the hypersphericity of $h_N$ in $T^N$ tending
to 0.  Taking surfaces $\Sigma_i$ in the class $h_N$ converging
to $T^N$, we get more examples of surfaces with Uryson 1-Width at
least 1 and hypersphericity tending to 0.

Arguments using maps of small girth and homotopy theoretic
obstructions cannot prove the small hypersphericity of these
surfaces $\Sigma_i$.  The reason is that, by the first extension
lemma, any map of small girth from $\Sigma_i$ to X extends to a
map from X to $T^N$ homotopic to the inclusion of $\Sigma_i$ in
$T^N$.  Since there is a continuous map from $T^N$ to $S^2$ not
killing the homology class of $\Sigma_i$, we see that there is a
continuous map from X to $S^2$ not killing the homology class of
$\Sigma_i$.  

Slightly stronger results are available using the Conway-Thomson
lattices.  Conway and Thomson proved that there are self-dual
lattices in $\mathbb{R}^N$ with shortest vector on the order of
$\sqrt{N}$.  (The theorem of Conway and Thomson is given as
Theorem 9.5 in \cite{MH}.)  If we take $T(N)$ to be the quotient
of $\mathbb{R}^N$ by a Conway-Thomson lattice, we get a torus
with systole on the order of $\sqrt{N}$, and with the property
that every homotopically non-trivial map to the unit circle has
Lipshitz constant at least on the order of $\sqrt{N}$.  This
property shows that the Conway-Thomson tori obey a stronger
version of Lemma 5.2:  the conclusion of the lemma holds for
every homology class h in $H_1(T(N))$.

Using the Conway-Thomson tori, we can produce metrics on high
genus surfaces with arbitrarily large homology systole and no
homotopically non-trivial contracting map to the unit circle.  In
particular, such a surface admits no degree 1 contracting map to
the unit square torus. (Recall that, by Theorem 4.1, any surface
with homology systole at least $4 \pi$ admits a degree 1
contracting map to the unit 2-sphere.)  We take a surface
$\Sigma$ of genus G embedded in a Conway-Thomson torus T(2G) as
its Jacobian torus.  The first homology group of $\Sigma$ is
identified with the first homology group of $T(2G)$.  We isotope
$\Sigma$ so that, for many homology classes $h_i$ in the first
homology group of $\Sigma$, there is a closed loop in the class
$h_i$ which lies very near to the straight circle in $T^N$ in the
corresponding homology class.  We put the induced metric on
$\Sigma$.  It is not hard to check that $\Sigma$ has (homology)
systole on the order of $\sqrt{N}$, and that every homotopically
non-trivial map from $\Sigma$ to the unit circle has Lipshitz
constant on the order of $\sqrt{N}$.  Taking N large proves our
claim.

The Conway-Thomson lattices were suggested to me by Mikhail Katz,
who was the first to apply them to systolic problems in his paper
\cite{Ka}.

We now turn to our last example of a manifold of bounded geometry
carrying a homology class of small hypersphericity.  This example
is based on some considerations arising in the theory of
systoles.  Let us define the length of an integral homology class
to be the smallest length of a union of curves realizing that
homology class.  On manifolds of dimension at least 3, it can
happen that the norm of a class $\alpha$ is arbitrarily large and
at the same time the norm of the class $100 \alpha$ is
arbitrarily small.  We give a simple construction of a Riemannian
manifold demonstrating this phenomenon.

Begin with the direct product $S^1 \times S^2$, with the standard
metric.  Let T be a small tubular neighborhood of an embedded
curve homologous to $100 [S^1]$.  On the interior of T, we modify
the metric by a conformal factor of $\epsilon$, and on the
exterior of T we modify the metric by a conformal factor of R,
where R is some enormous number depending on $\epsilon$.  (Near
the edge of T, we can smooth the jump so that the Riemannian
metric remains smooth.)  The length of $100 [S^1]$ is around
$200 \pi \epsilon$.  By taking R sufficiently large, we can make
the length of $[S^1]$ as large as we like.

This phenomenon can occur on manifolds with bounded geometry as
well.  We will construct a metric $g_N$ on $S^1 \times S^2$ with
bounded geometry on the scale 1, and with the property that the
homology class $N [S^1]$ has length only $2 \pi$.  We let T be a
small tubular neighborhood of an embedded curve homologous to $N
[S^1]$.  We put a Riemannian metric on $S^1 \times S^2$ which,
restricted to T, is just the product metric of the disk of radius
1 with the circle of radius 1.  Then we rescale this metric by a
large factor R, so that it has bounded geometry on the scale 1. 
Next, we will change the metric on the interior of T so that the
homology class $N [S^1]$ is realized by a curve of length $2
\pi$.  We choose coordinates on T so that the disk of radius R
has standard polar coordinates $r$ and $\theta$, and the circle
of radius R has a coordinate z varying from 0 to $2\pi$.  Our
rescaled metric on T is given in these coordinates as $g = dr^2 +
r^2 d\theta^2 + R^2 dz^2$.  We will replace this metric by a
metric $g_N = dr^2 + r^2 d\theta^2 + f(r)^2 dz^2$, where $f(r) = R$
for $R/2
\le r \le R$, and where $f(r) = 1$ for $0 \le r \le 1$.  By a
routine computation, it follows that the curvature of this metric
is bounded by 1 as long as the first and second derivatives of
log f are bounded by 1/2.  For large values of R, we can easily
find a (non-strictly) monotonic function f, satisfying the
equations above and with the first and second derivatives of log
f as small as we like.  Finally, since f is monotonic, any curve
in T contracts to a curve which runs along the core circle
defined by the equation r=0.  Since this core circle has length
$2 \pi$, we see that every closed geodesic in $(T, g_N)$ has
length at least $2 \pi$.  Given that the curvature of $g_n$ is
bounded by 1, the bound on the lengths of geodesics proves that
the part of $(T, g_N)$ a distance at least $\pi$ from the
boundary of T has injectivity radius at least $\pi$.  On the
other hand, $g_N$ has injectivity radius at least $\pi$ on the
rest of $S^1 \times S^2$, where it is equal to our original
rescaled metric.  Therefore, $g_N$ has bounded geometry on the
scale 1, while $N [S^1]$ is realized by a circle of length $2
\pi$.

Any homotopically non-trivial map from $(S^1 \times S^2, g_N)$ to
the unit circle must have Lipshitz constant at least N. 
Therefore, the hypersphericity of the class $[S^1]$ is at most
$1/N$.  If we take the product of this Riemannian manifold with
itself, we get a metric on $(S^1 \times S^1 \times S^2 \times
S^2)$ with bounded geometry on the scale 1, and the
hypersphericity of the homology class $[S^1 \times S^1]$ is
bounded by $\sqrt{\pi}/N$, because $(S^1 \times S^1 \times S^2
\times S^2, g_N \times g_N)$ contains a torus of area $4 \pi^2$
in the homology class $N^2 [S^1 \times S^1]$.  Any map to the
a sphere which does not kill $[S^1 \times S^1]$ induces a map
from this torus of degree at least $N^2$.  If the map is
contracting, then the area of the image sphere is bounded by $4
\pi^2 / N^2$, and so the hypersphericity of the class $[S^1
\times S^1]$ is bounded by $\sqrt{\pi}/N$.  As in previous
examples, we can construct a sequence of surfaces $\Sigma_i$ in
$(S^1 \times S^1 \times S^2 \times S^2)$, in the homology class
$[S^1 \times S^1]$, and converging in the Gromov-Hausdorff sense
to $(S^1 \times S^1 \times S^2 \times S^2, g_N \times g_N)$. 
Using Lemma 5.3, we see that these surfaces have Uryson 1-Width
at least 1 and arbitrarily small hypersphericity.

Finally, I would like to mention that the Lipshitz constant of a
non-zero degree map from $\Sigma_i$ to $CP^{\infty}$ is
arbritrarily large.  Unfortunately, Lemma 5.3 does not provide
such a bound, because the dimension of $CP^{\infty}$ is infinite. 
We can still prove the bound by an ad hoc trick.  First
triangulate the manifold $(S^1 \times S^1 \times S^2 \times S^2,
g_N \times g_N)$ into very small, approximately Euclidean simplices.  Say
that the simplices are all larger than some small number
$\epsilon$.  Then take a surface $\Sigma_i$ which is
Gromov-Hausdorff $\delta$-close to the total space, for some
$\delta$ much smaller than $\epsilon$.  Let f be a Lipshitz map from
$\Sigma_i$ to $CP^{\infty}$, with Lipshitz constant L.  By Lemma
5.2, we get a map from the 0-skeleton of the triangulation to
$CP^{\infty}$ with Lipshitz constant very close to 1.  We then
extend this map to all of $(S^1 \times S^1
\times S^2 \times S^2, g_N \times g_N)$ by straightening each simplex.  It
is not hard to check that we get a map with Lipshitz constant
less than $c L$, where c is a constant depending on the dimension
of our domain manifold (which is 6), but not depending on the
metric $g_N$ (because each extension occurs on an approximately
Euclidean simplex).  Suppose that f has degree non-zero.  Then,
we restrict the extension of f to the small area torus in the
homology class $N^2 [S^1 \times S^1]$.  The restriction is a
$cL$-Lipshitz map from a torus of area $4 \pi^2$ to $CP^{\infty}$
with degree at least $N^2$.  Therefore, L is greater than a
multiple of N, which proves our claim.

\end{document}